\documentclass{ectj}

\usepackage{amsmath,amsfonts,amssymb,graphics,epsfig,verbatim,bm,latexsym,url,amsbsy}
\usepackage{booktabs} 
\usepackage{siunitx}

\usepackage{amsthm}
\usepackage{xr}
\usepackage{mathabx}

\usepackage{xr}
\usepackage{natbib}
\usepackage{graphicx}
\usepackage{lscape}
\usepackage{subcaption}
\usepackage{color}
\usepackage{wrapfig}
\usepackage{mathtools}
\usepackage{placeins}
\usepackage{tgpagella}
\usepackage[framemethod=TikZ]{mdframed}
\usepackage{xcolor}
\usepackage{tikz}
\usetikzlibrary {positioning}

\usepackage{bm}
\usepackage{paralist}
\usepackage{accents}
\usepackage{enumerate}

\newtheorem{theorem}{Theorem}
\newtheorem{assumption}{Assumption}

\newtheorem{corollary}{Corollary}
\newtheorem{lemma}{Lemma}
\usepackage{tabularx}
\DeclareMathOperator{\eig}{eig}
\theoremstyle{remark}
\newtheorem{example}{Example}
\newtheorem{remark}{Remark}
\newtheorem{definition}{Definition}[section]

\renewcommand{\thesection}{\arabic{section}}
\renewcommand{\theequation}{\arabic{section}.\arabic{equation}}
\renewcommand{\thetheorem}{\arabic{section}.\arabic{theorem}}
\renewcommand{\theassumption}{\arabic{section}.\arabic{assumption}}
\renewcommand{\theproposition}{\arabic{section}.\arabic{proposition}}
\renewcommand{\thecorollary}{\arabic{section}.\arabic{corollary}}
\renewcommand{\thelemma}{\arabic{section}.\arabic{lemma}}
\renewcommand{\theexample}{\arabic{section}.\arabic{example}}
\renewcommand{\theremark}{\arabic{section}.\arabic{remark}}

\renewcommand{\Pr}{{\mathrm{P}}}

\newcommand{\mcL}{{\mathcal L}}

\newcommand{\Enk}{{\mathbb{E}_{n,k}}}
\newcommand{\Gnk}{{\mathbb{G}_{n,k}}}
\newcommand{\En}{{\mathbb{E}_n}}

\newcommand{\R}{\mathbb{R}}

\newcommand{\E}{\mathbb{E}}
\usepackage{nccmath}
		{\end{pmatrix}\end{medsize}}%

\renewcommand{\Pr}{{\mathrm{P}}}

\renewcommand{\hat}{\widehat}

\renewcommand{\Pr}{{\mathrm{P}}}

\renewcommand{\hat}{\widehat}
\renewcommand{\leq}{\leqslant}
\renewcommand{\geq}{\geqslant}

\newcommand{\mcH}{{\mathcal H}}

\newcommand{\mcP}{{\mathcal P}}

\newcommand{\mcG}{{\mathcal G}}
\newcommand{\mcF}{{\mathcal F}}

\newcommand{\mcZ}{{\mathcal Z}}
\newcommand{\mcW}{\mathcal W}

\received{February 2017}
\accepted{November 2017}
\volume{20}

\title[Regularized Machine Learning]{Regularized Orthogonal Machine Learning for Nonlinear Semiparametric Models}

\author[Nekipelov et al.]{Nekipelov, Denis$^{\dagger}$ and
                        Semenova, Vira$^{\ddagger}$ and Syrgkanis, Vasilis$^{\star}$\thanks{We are extremely thankful to the editor and three anonymous referees for the comments that helped us improve the paper. 
                        We thank           Victor Chernozhukov for many helpful discussions and pointing us to helpful references. We also thank  Patrick Kline, Whitney Newey, James Powell and seminar participants at Harvard/MIT econometrics seminar and Microsoft Research for useful comments.}}

\address{$^{\dagger}$ 254 Monroe Hall University of Virginia, Charlottesville, VA 22904, USA}
\email{denis@virginia.edu}

\address{$^{\ddagger}$530 Evans Hall, University of California, Berkeley, CA 94720}
\email{semenovavira@gmail.com}

\address{$^{\star}$ 1 Memorial Drive, Cambridge, MA, 02142 }
\email{vasy@microsoft.com}

\def\AmSTeX{$\cal A$\kern-.1667em\lower.5ex\hbox{$\cal M$}\kern-.125em
            $\cal S$-\TeX}
\def\BibTeX{{\rm B\kern-.05em{\sc i\kern-.025em b}\kern-.08em
            T\kern-.1667em\lower.7ex\hbox{E}\kern-.125emX}}

\begin{document}

\begin{abstract}
This paper proposes a Lasso-type estimator for a high-dimensional sparse parameter identified by a single index conditional moment restriction (CMR). In addition to this parameter, the moment function 
can also depend on a nuisance function, such as the propensity score or the conditional choice probability, which we estimate by modern machine learning tools.  We first adjust the moment function so that the gradient of the future loss function is insensitive (formally, Neyman-orthogonal)  with respect to the  first-stage regularization bias, preserving the single index property. We then take the  loss function to be an indefinite integral of the adjusted moment function with respect to the single index. The proposed Lasso estimator converges at the oracle rate, where the oracle knows the nuisance function and solves only the parametric problem. We demonstrate our method by estimating the short-term heterogeneous impact of Connecticut's Jobs First welfare reform experiment on women's welfare participation decision.

\end{abstract}

  \section{Introduction}
  \label{sec:intro}

Conditional moment restrictions (CMRs) often emerge as natural restrictions summarizing conditional independence, exclusion, or structural assumptions in economic models.  In such a setting, a major challenge is to allow for a data-driven selection among a (very) large number of conditioning covariates. The sparsity assumption, which requires  the number of relevant covariates to be small, has appeared to be an interpretable and plausible alternative. If the target parameter minimizes some population loss, a natural way to impose sparsity is to add an $\ell_1$-penalty on the target parameter to the sample loss. As discussed in \cite{EHT}, this Lasso approach has substantial computational and statistical advantages over its alternatives.  However, in general, the loss may be difficult to find. Focusing on a class of single index CMRs (\cite{Ichimura:93}, \cite{kleinspady}), we   make the Lasso approach feasible by deriving the  loss for an arbitrary single index moment function and establish convergence rate for the Lasso estimator.

The starting point of the analysis is a single index CMR that, in addition to target parameter, may depend on a nuisance component,  such as  the propensity score, the conditional choice probability, the conditional density, or alike, that can be substantially more high-dimensional than the target parameter itself. A natural approach would be to plug a machine learning/regularized estimate of the nuisance parameter into the moment function,  such as random forest or neural networks.  We first adjust the moment function so that the gradient of the future loss function is insensitive (i.e., Neyman-orthogonal, \cite{chernozhukov2016double}) with respect to the first-stage regularization bias, preserving the single index property.  
We then take the loss function to be an indefinite integral of the adjusted moment function with respect to the single-index. If the original moment function is monotone in the single-index, we report the global minimum of $\ell_1$-regularized $M$-estimator loss, following \cite{Negahban}.  Under mild conditions, the proposed estimator converges at the oracle rate, where the oracle knows the true value of the nuisance parameter and solves only the parametric problem.

We demonstrate the utility of our method with theoretical and  empirical applications.  First, we introduce a partially linear logistic model with heterogeneous treatment effects and derive an orthogonal loss. Second, we derive an orthogonal  loss  in conditional moment models with missing data, as studied in \cite{carroll:95}, \cite{carroll:91}, \cite{chen:08}, \cite{lee:95}, \cite{stepanski:93} and static games of incomplete information (e.g. see \cite{bajari:10} and \cite{bajari:13} among others). In all these settings, we give sufficient primitive  conditions on the nuisance parameters to achieve oracle convergence.  In the empirical application, we study the heterogeneous Jobs First effect on welfare participation decision via a partially linear logistic model. To detect treatment effect's heterogeneity, it is essential to use the orthogonal loss rather than the non-orthogonal one.

\paragraph{Literature review. }This paper is related to three lines of research: single index models, high-dimensional sparse models, and orthogonal/debiased inference based on machine learning methods. The first line of research concerns with the estimation of single index models (\cite{manski:75}, \cite{powell}, \cite{manski:85},  \cite{Ichimura:93}, \cite{kleinspady}, \cite{ichimura:hall}). This work focuses on a low-dimensional target parameter that can be treated as fixed. Focusing on smooth models, we allow the parameter's dimension to grow with sample size and even exceed it. We show that the single index property is sufficient to ensure the uniform convergence of the sample moments towards its population analog over an $\ell_1$-restricted ball, extending  the generalization bounds in the machine learning literature (see e.g., \cite{Shalev2014}) to single index CMRs.

The second line of research establishes the finite-sample bounds for a high-dimensional sparse parameter (\cite{belloni:11}, \cite{Negahban},  \cite{loh:13}, \cite{geer},  \cite{loh2017}, \cite{zhu:17}, \cite{Zhu2017}). Our contribution is to allow the loss function to depend on a functional nuisance parameter.  
In the convex case, we establish global convergence of the $\ell_1$-regularized $M$-estimator, following  \cite{Negahban}. One could follow a similar path to establish local convergence,  building on \cite{loh:13} and \cite{loh2017}, in the non-convex case. Focusing on the double robustness property,  \cite{tan2017regularized} and \cite{tan2018modelassisted} establish convergence guarantees in high-dimensional models that are  potentially misspecified. After we released the working paper version of this article (arxiv.ID 1702.06240), many methods have  proposed similar $M$-estimation approaches. Focusing on missing data, \cite{chakrabortty2019high} develops an $M$-estimator, relying on a  classical \cite{Robins} orthogonal score. We derive an orthogonal loss function for an arbitrary single index moment restriction, including  \cite{Robins} as a special case. The follow-up paper by \cite{foster2019orthogonal} extends our result to  $M$-estimators with  decomposable regularizers, including $\ell_1$ penalty as a special case. An alternative regularized minimum distance approach has been proposed in \cite{belloni2018highdimensional}.

The third line of research  obtains a $\sqrt{N}$-consistent and asymptotically normal estimator of a low-dimensional target parameter  $\theta$ in the presence of a  nonparametric nuisance parameter (\cite{Neyman:1959},  \cite{Neyman:1979}, \cite{HardleStoker1989},\cite{bickel:93}, \cite{NeweyStoker},  \cite{andrews1994},\cite{Newey1994}, \cite{Robins},  \cite{chen:03}). A statistical procedure is called Neyman-orthogonal if it is locally insensitive with respect to the estimation error of the first-stage nuisance parameter. Combining orthogonality and sample splitting, the Double Machine Learning
framework of  \cite{LRSP} and \cite{chernozhukov2016double}   has derived a root-N consistent asymptotically normal estimator of the target parameter based on the first-stage machine learning estimates. Extending this work, we establish convergence rates for $\ell_1$-regularized $M$-estimators whose  loss function gradient is an orthogonal moment. Next, we also contribute to the literature that derives  orthogonal moments starting from  non-orthogonal ones. Specifically, we construct a bias correction term for a nuisance parameter that is identified by a general conditional exogeneity restriction, covering  conditional expectation (\cite{Newey1994}) and conditional quantile (\cite{Newey:2018}) as leading special cases. An alternative approach based on automatic debiasing has been proposed in \cite{chernozhukov2021debiased} and \cite{chernozhukov2021automatic}.   Finally, we also contribute to the growing literature on orthogonal/doubly robust estimation based on machine learning methods (\cite{sasaki2018estimation}, \cite{chiang2019multiway}, \cite{sasaki2020unconditional}, \cite{Chiang}), in particular,  heterogeneous treatment effects estimation (\cite{nie2017}, \cite{CGST}, \cite{OprescuWu}, \cite{Lieli}, \cite{ZimLech}, \cite{Colangelo}, \cite{CherSem}). In contrast to this work, our main example features partialling out inside the argument of nonlinear link function, which, to the best of our knowledge, is  completely new.

This paper is organized as follows. Section \ref{sec:notation} introduces our main examples and gives a non-technical overview of the results. Section \ref{sec:moment} formally states our results. Section \ref{sec:applications} derives the concrete conditions for our applications. Section \ref{sec:emp} gives an empirical application.  Appendix \ref{sec:theory} generalizes Theorem \ref{cor:non-ortho-rate}  to the case of extremum estimators beyond $M$-estimators.  Appendix \ref{sec:proofs_moment} proves Theorem \ref{cor:non-ortho-rate} and Theorem \ref{lem:general}. Appendix  \ref{sec:ortho:proof} verifies the conditions of Theorem \ref{cor:non-ortho-rate} for each of the three applications.

 \section{Set-Up}
\label{sec:notation}

We start with the description of  a single index conditional moment restriction (CMR) framework. The CMR takes the form
 \begin{align}
    \E \bigg[m(W , \Lambda(Z, \gamma)'\theta_0, \gamma) \bigg|_{\gamma = g_0(Z)} \bigg| Z=z \bigg]~&=0, \label{eq:main}  \quad \forall z,
  \end{align} 
where  the first argument $W \in \mcW \subseteq \mathrm{R}^{\text{dim} W}$ is the data vector, the second argument $t \in \mathrm{R}$ is the single index, and the third one $\gamma \in \mathrm{R}^d$ is the output of the functional nuisance parameter. The  parameter of interest $\theta \in \mathrm{R}^p$ enters the moment function only via its inner product 
   \begin{align*}
t=\Lambda(Z,\gamma)'\theta,
\end{align*}
where the index function $\Lambda(z,\gamma): \mathcal{Z} \bigtimes \mathrm{R}^d \rightarrow \mathrm{R}^p$ is known up to $\gamma$. In many cases (e.g.,  \cite{Ichimura:93}, \cite{kleinspady}), the index function $\Lambda(z,\gamma)$ reduces to $\Lambda(z,\gamma)=z$.  Examples of the nuisance vector-function
 $$ g_0 = g_0(z) $$
 include the propensity score, the conditional choice probability, and the regression function, or a combination of these functions.  Given the CMR \eqref{eq:main}, our goal is to find a loss function $Q(\theta, \gamma): \mathrm{R}^p \bigtimes \mathrm{R}^d \rightarrow \mathrm{R}$  so that the true parameter value $\theta_0$ obeys
 \begin{equation}
    \label{eq:ld}
 \theta_0 = \arg \min_{\theta} Q(\theta, g_0).
\end{equation}

We find the loss function in two steps. We first solve the   ordinary  differential equation (ODE)
 \begin{align}
 \label{eq:ode}
 \dfrac{\partial }{\partial t} \ell (w, t, \gamma) = m(w, t, \gamma), \quad t  \in \R.
 \end{align}
We then plug $t=\Lambda(Z,\gamma)'\theta$ and $\gamma = g(Z)$ into the sketch of the loss function $\ell(w, t, \gamma)$ to obtain a population loss
\begin{align}
\label{eq:mestimatorfinal}
Q(\theta, g) = \E \bigg[ \ell (W, \Lambda(Z, g(Z))'\theta, g(Z)) \bigg],
\end{align}
whose gradient at $g=g_0$ and $\theta = \theta_0$ is 
 \begin{align}
\label{eq:mom} 
\nabla_{\theta} Q(\theta_0, g_0) = \E \bigg[m(W , \Lambda(Z, \gamma)'\theta_0, \gamma)\cdot \Lambda(Z, \gamma)  \bigg|_{\gamma = g_0(Z)}   \bigg].
\end{align} 
 In what follows, we refer to $\ell (w, t, \gamma)$ as the loss sketch, to $Q(\theta, g)$ and  $\nabla_{\theta} Q(\theta, g)$ as the population loss and the population gradient,    to
\begin{align}
\label{eq:sampleloss}
\widehat{Q}(\theta, \widehat{g}):=\dfrac{1}{n} \sum_{i=1}^n \ell (W_i, \Lambda(Z_i, \widehat{g}(Z_i))'\theta, \widehat{g} )
\end{align} and to $\nabla_{\theta} \widehat{Q}(\theta, \widehat{g})$ as the sample loss and the sample gradient.  The conditional independence assumption \eqref{eq:main} ensures that \eqref{eq:mom} is a valid moment equation for $\theta_0$. If $m(w, t, \gamma)$ is non-decreasing (non-increasing) in $t$, $\theta_0$ is the unique minimizer of $Q(\theta, g_0)$ ($-Q(\theta, g_0)$). 

The proposed  estimator $\widehat{\theta}$ has two stages. First, on the auxiliary sample, we construct an estimate $\widehat{g}$ of the nuisance parameter $g_0$, using a machine learning estimator capable of dealing with the high-dimensional covariate vector $Z$. Second, on the main sample, the target parameter's estimate is taken to be the minimizer of $\ell_1$-regularized sample loss.  For a fixed vector, its $\ell_2$ norm is denoted by $\| \cdot \|_2$, the $\ell_1$ norm is denoted by $\| \cdot \|_1$, the $\ell_{\infty}$ norm is denoted by $\| \cdot \|_{\infty}$, and $\ell_0$ norm is denoted by $\| \cdot \|_{0}$. 
%

\begin{definition}[Regularized $M$-Estimator]
	\label{def:alg:lasso}
 	Given  $(\widehat{g}(Z_i))_{i=1}^n$ and the penalty parameter $\lambda \geq 0$, define
	\begin{align}
	\label{eq:alg:lasso}
	\widehat{\theta} 		& =: \arg \min_{\theta \in \mathrm{R}^p} \widehat{Q}(\theta, \widehat{g}) + \lambda \| \theta \|_1.
	\end{align}
	\end{definition}

As discussed in \cite{Negahban}, the sample gradient $\nabla_{\theta} \widehat{Q} (\theta_0, \widehat{g})$ at $\theta = \theta_0$ summarizes the noise of the problem. If the population gradient \eqref{eq:mom} possesses the orthogonality property (\cite{Neyman:1959})
\begin{align}
&\dfrac{\partial}{\partial r} \nabla_{\theta} Q ( \theta_0, r (g - g_0) + g_0) = 0 \label{eq:orthogfinal}, \quad \forall r \in [0,1),
\end{align}
the biased first-stage estimation error $\widehat{g}(Z) - g_0(Z)$ has no first-order effect on the sample gradient. As a result, there exists a moderate penalty choice  $\lambda=\lambda_{\text{mod}}$  obeying
\begin{align}
\label{eq:moderate}
\lambda_{\text{mod}} = O \left( g_n^2 + \sqrt{\dfrac{\log p}{n}} \right)
\end{align}
that is sufficiently large to dominate the noise
\begin{align}
\label{eq:eventnoise}
\lambda/2 \geq \| \nabla_{\theta} \widehat{Q} (\theta_0, \widehat{g}) \|_{\infty}  \quad \text{ with probability } 1-o(1).
\end{align}
 If \eqref{eq:orthogfinal} does not hold,   the event \eqref{eq:eventnoise} requires an aggressive  choice $\lambda=\lambda_{\text{agg}}$
 obeying
\begin{align}
\label{eq:aggresive}
\lambda_{\text{agg}} = C \left( g_n + \sqrt{\dfrac{\log p}{n}} \right).
\end{align}
 To sum up, for $C$ sufficiently large, there exists an admissible penalty level $\lambda = \lambda_{\text{adm}} $ obeying \eqref{eq:eventnoise}
\begin{align}
\label{eq:admissible}
\lambda_{\text{adm}} = C \left(\sqrt{\dfrac{\log p}{n}} +B_0 g_n + g_n^2\right),
\end{align}
which reduces to $\lambda_{\text{mod}} $ if   \eqref{eq:orthogfinal} holds  (i.e., $B_0 = 0$) and to $\lambda_{\text{agg}} $ otherwise. For  an admissible penalty choice, Theorem \ref{cor:non-ortho-rate} establishes  the following bounds
\begin{align}
\label{eq:optimal}
\| \widehat{\theta} - \theta_0 \|_2 = O_P (\sqrt{k} \lambda) , \quad \| \widehat{\theta} - \theta_0 \|_1 = O_P  (k \lambda).
\end{align}
In particular, if $B_0 = 0$  and  $g_n^2= o( \log p/n)^{1/2}$ and $\lambda =\lambda_{\text{adm}}$ obeying \eqref{eq:admissible}, $\widehat{\theta}$ converges  at the oracle rate, where the oracle knows  the nuisance parameter $g_0$ and estimates only $\theta_0$.  In what follows, if the condition \eqref{eq:orthogfinal} holds, we refer to the   $m(w, t, \gamma)$ and  $\ell(w, t, \gamma)$ as the orthogonal moment and orthogonal loss, respectively.

We conclude this section by  studying the partially linear  logistic model. The model takes the form
  \begin{align}
\label{eq:nonlinearTE}
 \E \bigg[ Y -   G\left( D  \cdot ((1,X)'\theta_0) + f_0(X) \right)\mid D, X\bigg] =0,
 \end{align}
 where $D \in \R$ is a one-dimensional base treatment,  $X  \in \R^p$ is a vector of controls,  $Y$ is a binary outcome,  $W=(D,X,Y)$ is the data vector, and $G(t)$ is the logistic  link function. The treatment variable $D$ affects the outcome $Y$ via its interactions with the controls $(1,X)$. In addition, $X$ affects $Y$ via the confounding function $f_0(X)$ that enters \eqref{eq:nonlinearTE} in an additively separable way.  To make progress,  most papers (see e.g., \cite{belloni2016postselection}) require this function to be linear  so that $\theta_0$ and the nuisance parameter  can be  estimated under a joint sparsity assumption.  We describe below how to circumvent this bottleneck.

Inspired by  \cite{robinson:88}, we propose to partial out the controls inside the link function's argument
   \begin{align}
\label{eq:nonlinearTE2}
\E \bigg[ Y -   G\left( (D-p_0(X))  \cdot ((1,X)'\theta_0) + q_0(X) \right)\mid D, X\bigg] =0,
 \end{align}
 where
 \begin{align}
p_0(x) &= \E[ D \mid X=x] \label{eq:cexp}
\end{align}
is  the conditional expectation of the treatment and
\begin{align}
\label{eq:cpartindex}
q_0(x) = \E[ G^{-1} (\E [ Y | D, X]) \mid X=x]
\end{align}
is the conditional expectation of the link function's argument. In contrast to \eqref{eq:nonlinearTE}, the nuisance parameters $p_0(x)$ and $q_0(x)$ are identified separately from $\theta_0$ and permit a wider class of approaches to estimate them.  Thus, equation \eqref{eq:nonlinearTE2} is a special case of \eqref{eq:main} with  the conditioning vector $Z=(D,X)$, the moment function
\begin{align}
\label{eq:mpre}
m_{\text{pre}}(w, t, \gamma)= y - G( t +\gamma_2), 
\end{align}
the single index  $t= (d - \gamma_1) \cdot (1,x)'\theta$ and the nuisance parameter  $g_0(x) = \{ p_0(x),q_0(x) \}$ whose output is denoted by $\gamma = (\gamma_1, \gamma_2)$.

One may be tempted to proceed with the moment function \eqref{eq:mpre}.  Solving the ODE \eqref{eq:ode} gives
\begin{align}
\label{eq:lpre}
\ell_{\text{pre}}(w, t, \gamma)= - y \log G( t +\gamma_2) - (1-y) \log (1-G( t +\gamma_2)).
\end{align}
The sample loss  $\widehat{Q}(\theta, (\widehat{p}, \widehat{q}))$  in \eqref{eq:sampleloss} coincides with the  negative logistic likelihood. However, the population gradient \eqref{eq:mom} does not obey the orthogonality condition \eqref{eq:orthogfinal}
\begin{align}
&\dfrac{\partial}{\partial r} \nabla_{\theta}  Q ( \theta_0, r (q - q_0) + q_0) \nonumber \\
&= -\E G'((D-p_0(X)) ((1,X)'\theta_0) + q_0(X) ) \cdot (D-p_0(X))   (q (X) -q_0(X)) (1,X)'  \neq 0. \label{eq:nonorthogfinal}
\end{align}
As a result, the biased estimation error   $\widehat{q}(X)-q_0(X)$  has a first-order effect on  the sample gradient. Thus,  the   admissible penalty choice  reduces to $\lambda_{\text{adm}}=\lambda_{\text{agg}}$ in \eqref{eq:aggresive}, which makes the estimator's rate in \eqref{eq:optimal} a slow one.

To restore orthogonality, we reweigh the moment function as in \cite{belloni2016postselection}. The sketch of the new moment function is
\begin{align}
\label{eq:newsketch}
    m(w, t, \gamma):=\dfrac{y - G( t +\gamma_2)}{ \gamma_3 },
\end{align}
where $\gamma=(\gamma_1,\gamma_2,\gamma_3)$ corresponds to the output of $g_0(z) = \{ p_0(x),q_0(x),V_0(d,x)\}$ and the  weighting function $V_0(d,x)$ is the conditional variance
\begin{align}
\label{eq:vee1}
V_0(d,x) = G( (d-p_0(x)) \cdot ((1,x)'\theta_0) + q_0(x)) \cdot (1-G( (d-p_0(x)) \cdot ((1,x)'\theta_0) + q_0(x)))
\end{align}
 Solving   the ODE \eqref{eq:ode} gives
  \begin{align}
\label{eq:log:loss:orthog}
    \ell(w, t, \gamma) &= 
    -\dfrac{1}{\gamma_3} \big( y\cdot \log\left(G\left( t+ \gamma_2 \right)\right) + (1-y)\cdot \log\left(1-G((t+ \gamma_2))\right) \big).
\end{align}
The sample loss coincides with the negative weighted logistic likelihood (\cite{belloni2016postselection}). In contrast to \eqref{eq:nonorthogfinal}, the population gradient obeys the orthogonality condition \eqref{eq:orthogfinal}
\begin{align}
&\dfrac{\partial}{\partial r} \nabla_{\theta} Q ( \theta_0, r (q - q_0) + q_0) \nonumber \\
&= -\E \dfrac{G'((D-p_0(X)) ((1,X)'\theta_0 ) + q_0(X) )}{G'((D-p_0(X)) ((1,X)'\theta_0 ) + q_0(X) )} \cdot (D-p_0(X))  (q (X) -q_0(X)) (1,X)' = 0. \label{eq:nonorthogfinal2}
\end{align}
As a result,  the admissible penalty choice  \eqref{eq:admissible}  reduces to $\lambda_{\text{adm}}=\lambda_{\text{mod}}$ in \eqref{eq:moderate}, which makes the estimator's rate in \eqref{eq:optimal} a fast one.

\subsection{Examples}
\label{sec:examples}

\begin{example}[Nonlinear Treatment Effects]
\label{ex:nonlinearTE}

Suppose the treatment effect $\theta_0$ is identified by the conditional moment restriction \eqref{eq:nonlinearTE2}, where the link function $G(\cdot): \R \rightarrow \R$ is a known monotone link function that may not  necessarily be logistic. Define the weighting function as
\begin{align}
\label{eq:vee2}
V_0(d,x) = G'( (d-p_0(x)) \cdot ((1,x)'\theta_0) + q_0(x)).
\end{align}
 The  loss sketch $\ell(w,t,\gamma)$ is an arbitrary solution to the ODE
\begin{align}
\label{eq:nonlinearTEODE}
\dfrac{\partial}{\partial t}\ell(w,t,\gamma) = -\dfrac{1}{\gamma_3} \bigg(y-G( t +   \gamma_2  )  \bigg),
\end{align}
where  $\Lambda(z, \gamma) = (d - \gamma_1) \cdot (1,x)$,  $\gamma = (\gamma_1, \gamma_2, \gamma_3)$ denotes the output of $g_0(z) = \{ p_0(x), q_0(x), V_0(d,x)\}$ defined in \eqref{eq:cexp}, \eqref{eq:cpartindex} and \eqref{eq:vee2}. 
 Corollary \ref{cor:non-linear-te}  establishes the  mean square convergence rate   for the Regularized $M$-Estimator based on the loss sketch defined in  the ODE \eqref{eq:nonlinearTEODE}.

\end{example}

\begin{remark}[Linear Link]
\label{ex:linear-intro}
Consider Example \ref{ex:nonlinearTE} with $G(t)=t$.  The function $V_0(d,x)$ in \eqref{eq:vee2} simplifies to $$V_0(d,x)=1.$$ As a result, the nuisance parameter $g_0(z)$ simplifies to $g_0(z) = g_0(x) = \{ p_0(x), q_0(x)\}$ and $\gamma = (\gamma_1, \gamma_2)$.  The moment sketch is
\begin{align*}
    m(w,t, \gamma):= -(y - t-\gamma_2).
\end{align*}
The loss sketch  $\ell(w,t,\gamma)$ is 
\begin{align*}
\ell(w,t,\gamma) &= \dfrac{1}{2}(y  - t - \gamma_2)^2,
\end{align*}
which corresponds to the least squares loss used in \cite{CGST} and \cite{nie2017}. The population gradient  \eqref{eq:mom}  reduces to \cite{robinson:88}-type score 
\begin{align}
\label{eq:orthogfinal:linear}
  \nabla_{\theta} Q(\theta_0, g_0)  &= \E (Y  - (D-p_0(X)) \cdot (1,X)'\theta_0 - q_0(X)) \cdot  (D-p_0(X)) \cdot (1,X) =0.
\end{align}
Its pathwise derivative with respect to $q$ is zero:
\begin{align}
\label{eq:orthogfinal:linear2}
&\dfrac{\partial}{\partial r} \nabla_{\theta}  Q ( \theta_0, r (q - q_0) + q_0) =- \E  (D-p_0(X))   (q (X) -q_0(X)) \cdot (1,X)' =0. 
\end{align}
\end{remark}

\begin{remark}[Logistic Link]\label{ex:logistic}
Consider Example \ref{ex:nonlinearTE} with   $G(t)=(1+ \exp^{-t})^{-1}$.  The moment sketch is
\begin{align}
\label{eq:mlogistic}
    m(w, t, \gamma):=\dfrac{y - G( t +\gamma_2)}{ \gamma_3 }.
\end{align}
The loss sketch $\ell(w,t,\gamma)$ is  
\begin{align*}
    \ell(w, t, \gamma) &= 
    -\dfrac{1}{\gamma_3} \big( y\cdot \log\left(G\left( t+ \gamma_2 \right)\right) + (1-y)\cdot \log\left(1-G((t+ \gamma_2)\right) \big),
\end{align*}
which corresponds to the negative weighted logistic  likelihood used in \cite{belloni2016postselection}. The population gradient is 
\begin{align*}
   \nabla_{\theta} Q(\theta_0, g_0)   =-\E \dfrac{ (Y - G ((D-p_0(X)) \cdot (1,X)'\theta_0 + q_0(X)) }{V_0(D,X)} \cdot  (D-p_0(X)) \cdot (1,X),
\end{align*}
where $g_0(z) = (p_0(x),q_0(x), V_0(d,x))$ is as defined in \eqref{eq:cexp}, \eqref{eq:cpartindex}, \eqref{eq:vee1}.
\end{remark}

\begin{example}[Missing Data]
\label{ex:MD}

Suppose a researcher is interested in the parameter $\theta_0$ identified by a CMR:
\begin{equation}\label{eq:moment}
\E[u(Y^{*}, X'\theta_0) | X=x]=0, \quad \forall x \in \mathcal{X},
\end{equation}
where $Y^{*} \in \R$ is a partially observed outcome and $Z=X \in \R^p$ is a covariate vector. Let $V \in \{1,0\}$ indicate whether $Y^{*}$ is observed, $Y=V \cdot Y^{*}$ be the observed outcome, and $W=(V, X, Y)$ be the data vector.  A standard way to make progress is to assume that $V$ is as good as randomly assigned conditional on $X$.

Define the  conditional probability of observing $Y^{*}$ as
$$p_0(x) = \E[ V|X=x] $$
and the expectation function $h_0(x)$  as
 $$h_0(x) = \E\left[u(Y, X'\theta_0)\,|\,X=x, V=1\right].$$ The moment sketch is 
\begin{align}
\label{eq:orthoCMR2}
    m(w, t,  \gamma) = \dfrac{v\,}{\gamma_1}u(y, t)  -\dfrac{\gamma_2}{\gamma_1} ( v- \gamma_1),
\end{align}
where  $\gamma=(\gamma_1, \gamma_2)'$, $\Lambda(x, \gamma) = x$ and $g_0(x)=\{p_0(x), h_0(x)\}$.  The loss sketch is
\begin{equation}\label{eqn:m-loss-md}
\ell(w, t, \gamma) =   \dfrac{v\,}{\gamma_1}  \ell_{\text{pre}}(w, t,\gamma_1)  - \dfrac{\gamma_2}{\gamma_1}  (v-\gamma_1)\, t,
\end{equation} 
where $\ell_{\text{pre}} (w, t,\gamma_1)$ is an arbitrary solution to the ODE
\begin{align}
\label{eqn:m-loss-md:prelim}
\dfrac{\partial}{\partial t} \ell_{\text{pre}}(w,t,\gamma_1) =u(y, t).
\end{align}
The loss \eqref{eqn:m-loss-md} is a special case of the loss \eqref{eq:ortho} proposed in Theorem \ref{lem:general}. Corollary \ref{cor:selection}  establishes  mean square convergence rate   for the Regularized $M$-Estimator based on the loss function \eqref{eqn:m-loss-md}.

\end{example}

\begin{remark}[Quantile Regression with Missing Data]
\label{rm:quantile}
Consider Example \ref{ex:MD} with 
\begin{align}
\label{eq:quant}
u_{\tau}(y,t) = -(1_{ \{ y \leq t \}} - \tau),
\end{align} where $\tau \in (0,1)$ is a quantile level.  The moment function \eqref{eq:quant} identifies a \textit{quantile treatment effect} parameter $\theta_0$.  The loss sketch \eqref{eqn:m-loss-md} takes the form
 \begin{equation}\label{eqn:m-loss-md:quant}
\ell(w, t, \gamma) := - \left( \tau \cdot (y-t) 1_{ \{ y>t \}} + (1-\tau) \cdot (t-y) 1_{ \{ y<t \}} \right) \dfrac{v}{\gamma_1}  -\dfrac{\gamma_2}{\gamma_1} \,(v-\gamma_1)\, t.
\end{equation}

\end{remark}

\begin{example}[Static Games of Incomplete Information]
\label{ex:games}
Consider a two-player binary choice static game of incomplete information. The utility of action for player one is
\begin{align}
\label{eq:player1}
U (1) &= X'\alpha_0 + V \cdot \Delta_0+ \epsilon, \quad \E[ \epsilon | X] =0,
\end{align}
where  $X \in \R^p$ is the covariate vector, $V \in \{1, 0\}$ is the opponent action,  $\epsilon$ is mean independent private shock that follows Gumbel distribution.  The target $p$-vector $\theta_0=(\alpha_0, \Delta_0)$ consists of the covariate effect $\alpha_0$ and interaction effect $\Delta_0$.  The utility $U(0)$ of non-action for player one is normalized to zero. 

If the players' choices correspond to Bayes-Nash equilibrium, the outcome $Y$ obeys
\begin{align*}
Y &= 1[  X'\alpha_0 + p_0(X) \Delta_0 + \epsilon >0 ],
\end{align*}
where $p_0(x) = \E[ V|X=x]$. The moment sketch is
\begin{align*}
    m(w,t, \gamma) = -(y -G ( t ) + \gamma_2 (v - \gamma_1)),
\end{align*}
where  $\gamma=(\gamma_1, \gamma_2)'$, $\Lambda(x,\gamma) = (x, \gamma_1)$, and $g_0(x) = \{ p_0(x), h_0(x) \}$ for $h_0(x) = \Delta_0  G' ( x'\alpha_0 + p_0(x) \Delta_0 ).$   The loss sketch is
\begin{align}
\label{eq:losshardgames}
    \ell(w, t, \gamma ) :&=   \ell_{\text{pre}}(w, t, \gamma_1)-  \gamma_2\, (v - \gamma_1)\, t,
 \end{align}
 where the preliminary loss is the negative logistic likelihood
 \begin{align}
 \label{eq:losshardgames:pre}
 \ell_{\text{pre}}(w, t, \gamma_1) = - y  \cdot G (t) - (1-y) \cdot (1-G (t) ).
 \end{align} 
The loss \eqref{eq:losshardgames} is a special case of the loss \eqref{eq:ortho} proposed in Theorem \ref{lem:general}. 
 Corollary \ref{cor:games}  establishes mean square convergence rate for the Regularized $M$-Estimator based on the loss function \eqref{eq:losshardgames}. 

\end{example}

\subsection{ Remarks }
\label{sec:overview}


Following the sparsity bounds established in \cite{belloni2016postselection}, we conjecture that the final estimator $\widehat{\theta}$ is sparse. Therefore, one can interpret  the Lasso estimator $\widehat{\theta}$ as a model selector. We expect the post-Lasso-logistic based on single-selection procedure
\begin{align*}
 \widehat{\theta}_{PL} = \arg \min_{ \theta  \in \mathrm{R}^{p} }   \widehat{Q}(\theta, \widehat{g}): \quad  \text{support} (\theta) \subseteq  \text{support} (\widehat{\theta})
\end{align*}
to be too sensitive to moderate model selection mistakes, occurring when  the non-zero coefficients of $\theta_0$ are statistically indistinguishable from zero. According to \cite{LeebPotcher}, such inference procedures
do not provide a Gaussian approximation that is uniform over the space of $\theta$ and $g_0$ and are not honest. Instead, we discuss the following the debiasing procedure of \cite{geer}.

\begin{remark}[Debiased Lasso of \cite{geer}]
\label{rm:inference}
Suppose $t \rightarrow m(w, t, \gamma)$ is strictly increasing in $t$ for any $w$ and $\gamma$. Abstracting away from any nuisance components, or, effectively, treating $g_0$ as known, the work of \cite{geer} proposes a debiased Lasso estimator (eq. 18):
\begin{align}
\label{eq:debiased}
\widehat{\theta}_{\text{debiased}} (g_0)= \widehat{\theta} -\widehat{\Gamma} \frac{1}{n} \sum_{i=1}^n m (W_i, \Lambda(Z_i, g_0(Z_i))'\widehat{\theta},g_0) \cdot  \Lambda(Z_i, g_0(Z_i)),
\end{align} 
where  $\widehat{\theta} $ is a preliminary estimator of $\theta_0$, $\Gamma =(\nabla_{\theta \theta} Q(\theta, g_0))^{-1}$ is the population Hessian inverse, and $\widehat{\Gamma} $ is the estimator of $\Gamma$.  If the matrix $\Gamma$ is sparse, the estimator $\widehat{\Gamma} $ can be constructed by nodewise regression. Unlike the post-selection estimator,  $\widehat{\theta}_{\text{debiased}} (g_0)$ is Neyman-orthogonal with respect to the bias in the estimation error of $\widehat{\theta}-\theta_0$ and $ \widehat{\Gamma}-\Gamma_0$. If the sparsity indices of $\Gamma$ and the parameter $\theta_0$ are sufficiently small,   \cite{geer} shows that the estimator $\widehat{\theta}_{\text{debiased}} (g_0)$ is asymptotically Gaussian.  We conjecture that the plug-in estimator $$\widehat{\theta}_{\text{debiased}} (\widehat{g})$$ continues to be asymptotically Gaussian. For the linear link function, the asymptotic normality of $\widehat{\theta}_{\text{debiased}} (\widehat{g})$ has been established in \cite{CGST}. 
\end{remark}

\section{Theoretical results}
\label{sec:moment}

\paragraph{Notation. } We will use the following notation. For two sequences of random variables  $a_n, b_n, n \geq 1: a_n \lesssim_{P}  b_n$ means    $a_n = O_{P} (b_n)$. For two sequences of numbers $a_n, b_n, n \geq 1$, $a_n \lesssim  b_n$ means $a_n = O (b_n)$.  Let $T = \{ j: \quad \theta_{0, j} \neq 0 \}$ be the set of coordinates of $\theta_0$ that are not equal to zero, and let $T^c$ be the complement of $T$. For a vector $\delta \in \R^p$, let $(\delta_T)_j = \delta_j$ for each $j \in T$ and $(\delta_T)_j = 0$ for $j \in T^c$. 
For a vector-valued function $g(z):\mathcal{Z} \rightarrow \R^{d}$, denote its $\ell_2$ and $\ell_{\infty}$ norms as
\begin{align}
\label{eq:norms}
	 \| g \|_2 &:= (\E [ \| g(Z) \|_1]^2)^{1/2} =  (\E  (\sum_{l=1}^d | g_l(Z) | )^2)^{1/2} \\
	  \| g \|_{\infty} &:= \sup_{z \in \mathcal{Z}} \| g(z) - g(z') \|_1 = \sup_{z \in \mathcal{Z}}  \sum_{l=1}^d | g_l(z) - g_l(z') | \nonumber. 
\end{align} 
Furthermore, assume that the index function $\Lambda(z,\gamma): \mcZ \bigtimes \R^d \rightarrow \R^p$ is sufficiently smooth with respect to $\gamma$, so that the gradient $ \nabla_{\gamma}  \Lambda_j(z, \gamma)$  and the Hessian $ \nabla_{\gamma \gamma}  \Lambda_j(z, \gamma)$ of each coordinate $j \in \{1,2,\dots, p\}$ are well-defined. Finally, we will use the empirical process notation
 $$\En f(W_i) := \dfrac{1}{n} \sum_{i=1}^n f(W_i).$$ 

\begin{assumption}[Monotonicity in single index]
\label{ass:convexitymoment}
The moment function $m(w, t, \gamma)$ is non-decreasing in $t$ for any $w \in \mcW$ and any $\gamma \in \Gamma$.
\end{assumption}
Assumption  \ref{ass:convexitymoment} ensures that the loss sketch $\ell(w,t,\gamma)$ is a convex function of $t$ for any $w$ and $\gamma$. As a result, the sample loss $\theta \rightarrow \widehat{Q}(\theta, \widehat{g})$ defined in \eqref{eq:alg:lasso} is a convex function of $\theta$. By convexity,  on the event \eqref{eq:eventnoise},  the vector of errors
$
\nu = \widehat{\theta} - \theta_0 
$
belongs to the restricted cone
\begin{align}
\label{eq:cone}
 {\cal C}(T; 3) = \bigg\{ \nu \in \R^p: \quad \| \nu_{T^c} \|_1 \leq 3  \| \nu_{T} \|_1 \bigg\},
 \end{align}
 as shown in \cite{Negahban} (see  Lemma \ref{lem:subspace}). Define the restricted set as
\begin{align}
\label{eq:mathbb}
\mathbb{B} = \{\theta \in \R^p:  \theta = \theta_0 + r \nu, r \in [0,1), \nu \in {\cal C}(T;3) \}.
\end{align}
 
Define the sample curvature as 
\begin{align}
\label{eq:rsc}
\inf_{\theta \in \mathbb{B}, \nu=\theta - \theta_0 }  \dfrac{\nu^T  \nabla_{\theta \theta}  \widehat{Q}(\theta, \widehat{g})  \nu}{\| \nu \|_2^2},
\end{align}
and let the population curvature be the analog of \eqref{eq:rsc} based on the population Hessian $\nabla_{\theta \theta} Q(\theta, g_0) $ instead of $\nabla_{\theta \theta} \widehat{Q}(\theta, \widehat{g})$.

\begin{assumption}[Identification]
\label{ass:identification}
Let $\Sigma$ denote the population covariance matrix of the index function as
\begin{align}
\label{eq:sigma}
\Sigma &= \E \Lambda(Z, g_0(Z))\Lambda(Z, g_0(Z))^T.
\end{align}
Assume that  there exists a constant $C_{\text{min}}>0$ so that $\min \eig \Sigma \geq C_{\text{min}}$.
\end{assumption}
\begin{assumption}[Bounded derivative on $\mathbb{B}$]
\label{ass:identification2}
 There exists a constant $B_{\text{min}}>0$ so that the following bound holds: $$\inf_{\theta \in \mathbb{B} }\E \bigg[\nabla_{t}  m(W, t , \gamma)\bigg|_{\gamma = g_0(Z), t=\Lambda(Z, \gamma)'\theta}\,|\,Z=z \bigg] \geq B_{\text{min}}>0 \quad \forall z \in \mathcal{Z}. $$
\end{assumption}
Assumptions \ref{ass:identification} with $C_{\text{min}}$  and  \ref{ass:identification2} with $B_{\text{min}}$ imply that the population  curvature is bounded from below by
 $\bar{\gamma}  = B_{\text{min}} \cdot C_{\text{min}}$ (see Lemma \ref{lem:boundedhess}).  On the event
 \begin{align}
\label{eq:event}
\mathcal{V}:=\bigg\{  \sup_{\theta \in \mathbb{B}}  \| \nabla_{\theta \theta} [\widehat{Q}(\theta, \widehat{g} ) - Q(\theta,  g_0 )] \|_{\infty} < \bar{\gamma}/(32k)       \bigg\},
\end{align}
 the sample curvature \eqref{eq:rsc} is bounded from below by  $\bar{\gamma}/2$ (see Lemma \ref{lem:fsrsc}). Since  
  \begin{align*}
\nabla_{\theta \theta} [\widehat{Q}(\theta, \widehat{g} ) - Q(\theta,  g_0 )]  &\leq  \| \nabla_{\theta \theta} [\widehat{Q}(\theta, g_0 ) - Q(\theta,  g_0 )] \|_{\infty} + \| \nabla_{\theta \theta} [\widehat{Q}(\theta, \widehat{g} ) - \widehat{Q}(\theta, g_0 ) ] \|_{\infty},
\end{align*}
it suffices to show that each summand above is $o_P(1)$. We bound the first and the second summand in Lemmas \ref{lem:grc} and \ref{lem:liphessian-non-conv}, respectively.

 Assumption \ref{ass:smooth} ensures that the moment function is sufficiently smooth in its second and third arguments.  Let $\mcW$ be an open bounded set containing the support of the data vector $W$. Likewise, let $\mathcal{T}$ be an open set containing the support of  $\Lambda(Z,g_0(Z))'\theta_0$.   Finally, let $\Gamma$  be an open bounded set containing the support of vector $g(Z)$, when $g \in \mcG_n$. In what follows, we assume that the sets $\mcW, \mathcal{T}, \Gamma$ do not change with $n$.
  
\begin{assumption}[Smooth and bounded design]\label{ass:smooth}

There exists a constant $U< \infty$ so that $\sup_{\theta \in \Theta} \| \theta \|_1 \leq U$ and for any vector $w\in\,\mcW$, number $\ t\in\,\mathcal{T}$ and vector $\gamma \in \Gamma$ the following conditions hold:
    \begin{align*}
         |m(w, t, \gamma )|, \|\nabla_{\gamma} m(w, t, \gamma )\|_{\infty},\,   |\nabla_{t} m(w, t, \gamma)| \leq~& U,\\
         \|\nabla_{\gamma t} m(w, t, \gamma )\|_{\infty},     \| \nabla_{\gamma \gamma}  m(w, t, \gamma )\|_{\infty}     , |\nabla_{tt} m(w, t, \gamma)|\leq~& U,\\
             \|\Lambda_j(z, \gamma)\|_{\infty}, \|\nabla_{\gamma} \Lambda_j(z, \gamma)\|_{\infty}, \|\nabla_{\gamma \gamma}  \Lambda_j(z, \gamma)\|_{\infty} \leq~& U, \quad \forall j \in \{1,2,\dots, p\}.
    \end{align*}

\end{assumption}

Lemma \ref{lem:grc} establishes a convergence rate for the sample Hessian uniformly over the restricted set \eqref{eq:cone}. A gradient version of Lemma \ref{lem:grc} is available in the concurrent work by \cite{belloni2018highdimensional}.

\begin{lemma}[Uniform Convergence of Sample Hessian]
\label{lem:grc}
Suppose Assumption \ref{ass:smooth} holds with a constant $U$. Then,  there exists a sequence $\tau_n = 4 k U^5   \sqrt{\dfrac{2\log p}{n}} + 2U^3  \sqrt{\dfrac{2\log p}{n}}=o(1)$ so that $\nabla_{\theta \theta} \widehat{Q}(\theta, g_0) $ uniformly converges  to $\nabla_{\theta \theta} Q(\theta,g_0)$:
\begin{align}
\sup_{\theta \in \mathrm{B}} \| \nabla_{\theta \theta} \widehat{Q}(\theta, g_0) - \nabla_{\theta \theta} Q(\theta, g_0) \|_{\infty} = O_{P} (\tau_n) = o_P(1).
\end{align}
\end{lemma}

\begin{proof}
Consider the function class $$\mcF=\{Z \rightarrow \Lambda(Z, g_0(Z))'\theta, \quad \theta \in \R^p,  \quad \|\theta\|_1\leq k U  \}.$$ Define its Rademacher complexity 
$$
\mathcal{R}(\mcF) := \E \bigg[ \sup_{\theta: \|\theta\|_1\leq k U} \En \sigma_i \Lambda(Z_i, g_0(Z_i))'\theta \bigg], 
$$
where $(\sigma_i)_{i=1}^{n}$ is a vector of i.i.d random variables from Rademacher  distribution: $\Pr (\sigma_i=+1) = \Pr (\sigma_i=-1)=0.5$. Each function in the class   $$\mcH = \{ W \rightarrow \nabla_{\theta_k \theta_j} \ell (W, \Lambda(Z, g_0(Z))'\theta, g_0(Z)), \quad \theta \in \R^p,  \|\theta\|_1\leq k U \}$$ 
is a combination of the function in the class $\mcF$ and a $U^3$-Lipshitz function. Invoking  Contraction Lemma (Lemma 26.9) and Lemma 26.11 from \cite{Shalev2014}, we obtain
$$
\mathcal{R}(\mcH) \leq U^3 \mathcal{R}(\mcF) \leq  U^5 \cdot  k \sqrt{\frac{2\log(2p)}{n}}.
$$
Finally, invoking  Lemma 26.5 \cite{Shalev2014}  for each pair $(k, j) \in \{1,\ldots, p\}^2$, w.p. $1-\delta/p^2$,
\begin{align*}
    &\sup_{ \|\theta\|_1\leq k U} \bigg| [\En - \E ][\nabla_{\theta_k \theta_j} \ell(W_i, \Lambda(Z_i, g_0(Z_i) )' \theta_0, g_0(Z_i)) ] \bigg| \\
    &\leq 2 \mathcal{R}(\mcH) + U^3 \sqrt{\frac{2\log(2\,p^3/\delta)}{n}}\\
     &\leq 2   k U^5 \sqrt{\frac{2\log(2p^3/\delta)}{n}} + U^3 \sqrt{\frac{2\log(2\,p^3/\delta)}{n}}.
\end{align*}
Union bound over $p^2$ pairs $(k, j) \in \{1,\ldots, p\}^2$ and $2 p^3 \leq p^4$ gives w.p. $1-\delta$
\begin{align*}
    &\sup_{ 1 \leq k,j \leq p} \sup_{ \|\theta\|_1\leq k U} \bigg| [\En - \E ][\nabla_{\theta_k \theta_j} \ell(W_i, \Lambda(Z_i, g_0(Z_i) )' \theta_0, g_0(Z_i)) ] \bigg| \\
     &\leq 4 k U^5 \sqrt{\frac{2\log(p/\delta)}{n}} + 2U^3 \sqrt{\frac{2\log(p/\delta)}{n}} \
\end{align*}

\end{proof}

Assumption \ref{ass:nuisance} formalizes the convergence of the nuisance parameter's estimator. It introduces a sequence of nuisance realization sets $\mcG_{n} \subseteq \mathcal{G}$ that contain the true value $g_0$ and the  estimator $\widehat{g}$ with probability $1-\delta_n$. As the  sample size $n$ increases, the sets $\mcG_{n}$ shrink. The shrinkage speed is measured by the rate $g_{n}$ and is referred  to as the first-stage rate.

\begin{assumption}[First-stage rate]
\label{ass:nuisance} 
There exist sequences of numbers $\delta_n = o(1)$ and $g_n = o(1)$ and a sequence of sets $\mcG_n \subseteq \mcG$ so that $\widehat{g} \in \mcG_n$ w.p. at least $1-\delta_n$ and $g_0 \in \mcG_n$.   The sets shrink at the rate 
\begin{align*} 
\sup_{g \in \mcG_n} \| g - g_0 \| \leq g_{n},
\end{align*}
where $\| \cdot \|$ is either the $\ell_{\infty}$ norm or the $\ell_2$ norm as defined in equation \eqref{eq:norms}.
\end{assumption}

Assumption \ref{ass:nuisance} is satisfied by many machine learning estimators under structural assumptions on the model in $\ell_2$ and/or $\ell_{\infty}$ norm.  For example, it holds for Lasso (\cite{belloni:11}, \cite{belloni2013}) in linear and generalized linear models, for $L_2$-boosting in sparse models (\cite{Luo}), for neural network (\cite{ChenWhite1991}), and for random forest in low-dimensional (\cite{wagerwalther}) and high-dimensional sparse (\cite{fastRF}) models. For a broad overview of low-level primitive conditions that covers neural networks and random forest see, e.g., \cite{jeong2020robust}, Appendix 1.

\begin{theorem}[Regularized $M$-Estimator]\label{cor:non-ortho-rate}
Suppose that the nuisance space $\mathcal G$ is equipped with either the $\ell_{\infty}$ or
the $\ell_2$ norm, defined in equation \eqref{eq:norms}. Suppose Assumptions \ref{ass:convexitymoment}-\ref{ass:nuisance}  hold, and $k \left(g_n + U^5 \cdot k \cdot \sqrt{\dfrac{ 2\log 2p}{n}} + 3  d \cdot U^3 \, \sqrt{\frac{\log(2p)}{n}}\right)=o(1)$. Then, for $n$ and $C$ large enough and the admissible penalty $\lambda = \lambda_{\text{adm}}$ is chosen to obey \eqref{eq:admissible}, Regularized $M$-Estimator  obeys the bound  \eqref{eq:optimal}.

\end{theorem}

Theorem \ref{cor:non-ortho-rate} is our main result. It establishes the convergence rate for the Regularized $M$-estimator, covering the orthogonal (i.e., $B_0 = 0$) and the non-orthogonal (i.e., $B_0 \neq 0$)  cases. In the former case, the admissible penalty choice   $ \lambda_{\text{adm}}$ reduces to  $ \lambda_{\text{mod}}$ in \eqref{eq:moderate}, which makes the estimator's rate \eqref{eq:optimal} fast. Otherwise,  $ \lambda_{\text{adm}}$ reduces to  $ \lambda_{\text{agg}}$ in \eqref{eq:aggresive}, which makes the estimator's rate \eqref{eq:optimal} slow.  In all our examples, we invoke  Theorem \ref{cor:non-ortho-rate} twice: first, for  the preliminary estimate $\check{\theta}$ based on a non-orthogonal CMR \eqref{eq:main:pre} and, second, for the final estimate $\widehat{\theta}$ based on the orthogonal CMR.

Suppose that the moment function $m_{\text{pre}} (w, t, \gamma_1): \mcW \bigtimes \R \bigtimes \R \rightarrow \R$ and/or the index function $\Lambda(z, \gamma_1): \mcW \bigtimes \R  \rightarrow \R^p$ depend on a one-dimensional functional nuisance parameter $p_0(z)$:
\begin{align}
\label{eq:main:pre}
\E \bigg[ m_{\text{pre}} (W, \Lambda(Z, p_0(Z))'\theta_0, p_0(Z)) \bigg| Z=z \bigg] = 0, \quad z \in \mcZ. 
\end{align}
Furthermore, suppose the nuisance parameter is  identified by a  conditional exogeneity restriction
\begin{align}
\label{eq:cond_exog}
\E \left[ R(W, p_0(Z)) | Z=z \right] =0.
\end{align}
For example, $ R(W, p(X)) = D - p(X) $ defines the conditional expectation function as in \eqref{eq:cexp}, and $ R(W,p(X)) = 1_{ \{ V \leq p(X) \}} - \tau$ defines the conditional $\tau$-quantile function. Starting  from an arbitrary CMR \eqref{eq:main:pre}, Theorem \ref{lem:general} derives a loss whose gradient   obeys the orthogonality condition \eqref{eq:orthogfinal}.

\begin{theorem}[Construction of Orthogonal Loss]

\label{lem:general}
Suppose  equation \eqref{eq:cond_exog} holds. Define the adjusted moment function
 \begin{equation}
     \label{eq:rho:phi}
    m(w, t, \gamma) = m_{\text{pre}}(w , t, \gamma_1) - \gamma_2 \cdot \gamma_3^{-1} \cdot R(w, \gamma_1),
 \end{equation}
 where $\gamma = (\gamma_1, \gamma_2, \gamma_3)$ denotes the output of the nuisance parameter  $g_0(z) = \{ p_0(z), h_0(z), I_0(z) \}$, consisting of $p_0(z)$  as defined in \eqref{eq:cond_exog},  $h_0(z)$ and $I_0(z)$ defined as
  \begin{align}
 h_0(z)&=\E\bigg[ \nabla_{\gamma_1} m (W,t,\gamma_1) \bigg|_{\gamma_1=p_0(Z), t= \Lambda(Z, \gamma_1)'\theta_0} \bigg|Z=z\bigg] \\
 I_0(z) &= \E \bigg[ \nabla_{\gamma_1} R(W,\gamma_1)  \bigg|_{\gamma_1 = p_0(Z)} \bigg|Z=z \bigg].
 \end{align}
The  loss sketch $\ell(w,t, \gamma)$ takes the form
   \begin{align}
 \label{eq:ortho}
 \ell(w, t, \gamma) &= \ell_{\text{pre}}(w,t,\gamma_1) -  \gamma_2 \cdot \gamma_3^{-1} \cdot R(w, \gamma_1)  \cdot t,
 \end{align}
 where $\ell_{\text{pre}}(w,t,\gamma_1)$ solves the Ordinary Differential Equation \eqref{eq:ode} for the moment function $m_{\text{pre}} (w,t,\gamma_1)$. Then, the population loss $Q(\theta, g_0)$ defined in \eqref{eq:ld} obeys the orthogonality condition \eqref{eq:orthogfinal}. 
 
 Suppose  $\sup_{w \in \mcW} \sup_{\gamma_1 \in \Gamma_1} |R (w, \gamma_1)| \leq U$,  $\sup_{w \in \mcW} \sup_{\gamma_1 \in \Gamma_1} |\nabla_{\gamma_1} R (w, \gamma_1)| \leq U,$ $$\sup_{w \in \mcW} \sup_{\gamma_1 \in \Gamma_1}  |\nabla_{\gamma_1 \gamma_1} R (w, \gamma_1)| \leq U.$$   If the original moment function $m_{\text{pre}}(w, t, \gamma_1)$ obeys Assumptions \ref{ass:convexitymoment}-\ref{ass:smooth}, the adjusted moment function \eqref{eq:rho:phi} obeys Assumptions \ref{ass:convexitymoment}-\ref{ass:smooth}.  
 \end{theorem}

Starting from an arbitrary CMR \eqref{eq:main:pre},  Theorem \ref{lem:general} adjusts the  moment function so that population gradient $\nabla_{\theta} Q(\theta_0, g_0)$ obeys orthogonality condition \eqref{eq:orthogfinal}.  Special cases of equation \eqref{eq:cond_exog}, such as the  conditional expectation function   and conditional quantile function,   are available in \cite{Newey1994}, \cite{LRSP} and  \cite{Newey:2018} for unconditional moment problems. We extend the results above to allow  for  an arbitrary conditional exogeneity restriction.

\section{Applications}
\label{sec:applications}

In all the applications below, our starting point is the CMR \eqref{eq:main:pre} which does not obey the orthogonality condition \eqref{eq:orthogfinal}. Invoking either the weighting idea (Example \ref{ex:nonlinearTE}) or Theorem \ref{lem:general} (Examples \ref{ex:MD}-\ref{ex:games}), we derive an orthogonal CMR obeying \eqref{eq:orthogfinal}. However, this CMR requires an additional preliminary estimate of $\theta_0$ on top of the nuisance parameters involved in \eqref{eq:main:pre}. Definition \ref{def:twostep} describes the suggested two-step procedure, spelling out the choices of main/auxiliary samples  and the penalty parameters in each step.

\begin{definition}[Two-Step Regularized $M$-estimator]
\label{def:twostep}

Let $K=3$ denote the $3$-fold partition of the sample $\{1,2,\dots, n\}$ into $J_1, J_2, J_3$. For notational convenience, let $J_0 = J_3$ and $J_4=J_1$.  For each $k \in \{0,1,2\}$, compute
\begin{enumerate}
\item With $J_{k-1}$ as the auxiliary sample and $J_k$ as the main sample, let $\check{\theta}$ be  the output of Definition \ref{def:alg:lasso} with the preliminary loss $\ell_{\text{pre}}(w, t, \gamma_1)$ and the penalty parameter $\lambda_{\text{agg}}$ as in \eqref{eq:aggresive}. For each  $i \in J_{k+1}$, the nuisance parameter $\widehat{g}(Z_i)$ is evaluated at $\widehat{g}(Z_i) =\widehat{g}_k(Z_i) = \{  \widehat{p}_{k-1}(Z_i), \check{\theta}_k \}$.

\item  With ($J_{k-1}, J_k)$ as the auxiliary sample and $J_{k+1}$ as the main sample, let $\widehat{\theta}$ be the output of Definition \ref{def:alg:lasso} with the final loss $\ell(w, t, \gamma)$ and the admissible penalty parameter $\lambda_{\text{mod}}$ obeying \eqref{eq:moderate}. Report: $\widehat{\theta}$

\end{enumerate}

\end{definition}

\subsection{Nonlinear Treatment Effects}
\label{sec:nte1}
Consider  Example \ref{ex:nonlinearTE}.  Define the covariance matrix 
  \begin{equation}
  \label{eq:sigmate}
    \Sigma_{\text{TE}} := \E[(D-p_0(X))\, (D-p_0(X))(1,X)(1,X)^T].
     \end{equation}

 \begin{assumption}[Regularity Conditions for Nonlinear Treatment Effects]

 \label{ass:reg:nonlinearte}
 Suppose the following conditions hold.
\begin{enumerate}
       \item  \label{itm:nonlinear:ide} (Identification) There exists  a constant $\gamma_{\text{TE}}>0$ so that $\min \eig  \Sigma_{\text{TE}}  \geq \gamma_{\text{TE}}$.
       
        \item \label{itm:nonlinear:smooth} (Smooth and Bounded Design). The parameter space $ \Theta$ is bounded in $\ell_1$ norm by a constant  $H_{\text{TE}} \geq 1$:  $\sup_{\theta \in \Theta} \| \theta \|_1 \leq H_{\text{TE}}$ and $\sup_{w \in \mcW} \| w \|_{\infty} \leq H_{\text{TE}}$.   The functions $t \rightarrow G(t), G'(t), G^{(2)}(t), G^{3}(t)$ are $U_{\text{TE}}$-bounded on the set $[- 3 \cdot H_{\text{TE}}^3,3 \cdot H_{\text{TE}}^3]$, where $U_{\text{TE}} \geq  H_{\text{TE}}$.       Furthermore, the functions $t \rightarrow G^{-1}(t), (G')^{-1}(t), (G')^{-2}(t), (G')^{-3}(t)$ are $U_{\text{TE}}$-bounded  from above on the set $[- 3 \cdot H_{\text{TE}}^3, 3 \cdot H_{\text{TE}}^3]$.   Finally, $t \rightarrow G'(t)$ is $\text{B}_{\text{min}}$-bounded from below on the set  $[- 3 \cdot H_{\text{TE}}^3, 3 \cdot H_{\text{TE}}^3]$.

    \item \label{itm:nonlinear:rate} (First-Stage Rate). There exists an estimator $(\widehat{p}(x), \widehat{q}(x))$ of $(p_0(x), q_0(x))$ satisfying  Assumption \ref{ass:nuisance} with $\pi_{n}+q_{n}$ rates in either $\ell_{\infty}$ or $\ell_2$  norm.
    \item \label{itm:nonlinear:mon}  (Monotonicity). The function $G(t): \R \rightarrow \R$ is a known monotone function of $t$.
\end{enumerate}
\end{assumption}

\begin{remark}[Weighting function estimator]
\label{rm:fsvee}
For $k \in \{0,1,2\}$ and  $i \in J_{k+1}$,  define the estimator  $\widehat{V}(D_i,Z_i )=\widehat{V}_k(D_i,Z_i)$ 
\begin{align}
\label{eq:vee}
\widehat{V}_k(D_i,X_i)=G'((D_i-\widehat{p}_k(X_i)) \cdot (1,X_i)'\check{\theta}_k + \widehat{q}_k(X_i)), \quad i \in J_{k+1}.
\end{align}
where preliminary loss $\ell_{\text{pre}}(w, t, \gamma_1)$ as in \eqref{eq:lpre} is used for $\check{\theta}$.
\end{remark}

\begin{corollary}[Nonlinear Treatment Effects]
\label{cor:non-linear-te}
Suppose  Assumption  \ref{ass:reg:nonlinearte} holds. Let the preliminary loss $\ell_{\text{pre}}(w, t, \gamma_1)$ be as in \eqref{eq:lpre} and the first-stage parameter estimate $(\widehat{\pi}(x), \widehat{q}(x))$. Let the final loss be as in  \eqref{eq:nonlinearTEODE} and $\widehat{g} (d,x)= (\widehat{\pi}(x), \widehat{q}(x), \widehat{V}(d,x))$, where $\widehat{V}(d,x)$ as defined in  \eqref{eq:vee}.Then, the statement of Theorem \ref{cor:non-ortho-rate} holds for each step of  Definition \ref{def:twostep}, where the first-step nuisance rate is $O(\pi_{n}+q_{n})$ and the second-step nuisance rate is $g_n =  O \left( k  \left( \pi_n + q_n+\sqrt{\frac{\log p}{n}} \right) \right)$.
\end{corollary}

\subsection{Missing Data}
Consider  Example \ref{ex:MD}. Define the covariance matrix 
\begin{align}
\label{eq:sigmamd}
    \Sigma_{MD}= \E X X^T
    \end{align}
 and the function $q(t, x)$ as
 \begin{align*}
 q_0(x,t) = \E [ \nabla_t u( Y, t ) | V=1,  X=x]. 
 \end{align*}

\begin{assumption}[Regularity Conditions for Missing Data]
\label{ass:selection}
Suppose the following conditions hold. 
\begin{enumerate}
    \item \label{index-sel:identification} (Identification). There exists $\gamma_{\text{MD}}>0$  such that $\min \eig      \Sigma_{MD} \geq \gamma_{\text{MD}}$. Furthermore, there exists $B_{\text{min}}>0$ so that $\inf_{ t \in \mathcal{T}} \E [ \nabla_{t} u(Y,t) | X=x] \geq B_{\text{min}}$ for any $x \in \mathcal{X}$. 
    
       \item \label{index-sel:smooth} (Overlap Condition). There exists a constant $\bar{p}>0$ such that $\inf_{x \in \mathcal{X} } p_0(x)\geq \underline{p}>0$.

     \item  \label{index-sel:smooth}  (Smooth and Bounded Design).  The parameter space $ \Theta$ is bounded in $\ell_1$ norm by a constant  $H_{\text{MD}}$:  $\sup_{\theta \in \Theta} \| \theta \|_{\infty} \leq H_{\text{MD}}$ and $\sup_{w \in \mcW} \| w \|_{\infty} \leq H_{\text{MD}}$. The functions $u(y,t)$,$\nabla_{t} u(y,t)$, $\nabla_{tt} u(y,t)$ are $U_{\text{MD}}$-bounded in $\ell_{\infty}$ norm for all $w \in \mcW$ and $t \in \mathcal{T}$, where  $U_{\text{MD}}\geq H_{\text{MD}}$.
     
    \item \label{index-sel:nuisance} (First-Stage Rate). There exists an estimator $(\widehat{p}(x), \widehat{q}(x,t))$ of $(p_0(x), q_0(x,t))$ satisfying  Assumption \ref{ass:nuisance} with $\pi_{n,r}+q_{n,r}$ rates in either $\ell_{\infty}$ or $\ell_2$  norm, where $q_{n,r}$ rate is
    \begin{align*}
    q_{n, \infty}:= \sup_{q \in \mathcal{Q}_n} \sup_{t \in R} \sup_{x \in \mathcal{X}}  |q(x,t) - q_0(x,t)|, \quad  q_{n, 2}:=\sup_{q \in \mathcal{Q}_n}\sup_{t \in R} (\E (q(X,t) - q_0(X,t))^2)^{1/2}.
    \end{align*}

    \item \label{index-sel:convexity} (Monotonicity). The function $u(y,t)$ is increasing in $t$ for any $y \in \R$.

\end{enumerate}

\end{assumption}

\begin{remark}[First-Stage Estimator of $h_0(x)$]
\label{rm:fsvee:selection}
For each $i \in J_{k+1}$,  define the estimator $\widehat{h}(X_i)=\widehat{h}_k(X_i)$ 
\begin{align}
\label{eq:vee:selection}
\widehat{h}(X_i)=\widehat{q}_k(X_i'\check{\theta}_k, X_i),
\end{align}
where the preliminary estimator $\check{\theta}$ is based on the loss $\ell_{\text{pre}}(w, t, \gamma_1)$ as in \eqref{eqn:m-loss-md:prelim}.

\end{remark}

\begin{corollary}[General Moment Problems with Missing Data]
\label{cor:selection}

Suppose  Assumption  \ref{ass:selection} holds. Let the preliminary loss $\ell_{\text{pre}}(w, t, \gamma_1)$ be as in \eqref{eqn:m-loss-md:prelim}  and the first-stage parameter estimate $\widehat{p}(x)$. Let the final loss be as in  \eqref{eqn:m-loss-md}  and $\widehat{g} (x)= (\widehat{p}(x), \widehat{h}(x))$, where $\widehat{h}(x)$ as defined in  \eqref{eq:vee:selection}.  Then, the statement of Theorem \ref{cor:non-ortho-rate} holds for each step of  Definition \ref{def:twostep}, where the first-step nuisance rate is $O(\pi_{n})$ and the second-step nuisance rate $g_n =  O \left(k\left( \pi_n +\sqrt{\frac{\log p}{n}} \right)+q_n\right)$.

\end{corollary}

\subsection{Static Games of Incomplete Information}
\label{sec:strategic}
\label{SEC:STRATEGIC}
Consider  Example \ref{ex:games}. Define the covariance matrix 
    \begin{align}
       \Sigma_{\text{games}} &= \E[ (X,p_0(X)) (X, p_0(X))^T].
    \end{align}

\begin{assumption}[Regularity Conditions for Static Games of Incomplete Information]
\label{ass:games}
Suppose the following conditions hold.
\begin{enumerate}
    \item \label{index-games:identification} (Identification). There exists $\gamma_{\text{games}}>0$ so that $\min \eig     \Sigma_{\text{games}}  \geq \gamma_{\text{games}}$.
    
    \item \label{index-bounded} (Smooth and Bounded Design). There exists a constant $H_{\text{games}} < \infty$ so that  $\sup_{\theta \in \Theta} \| \theta \|_{\infty} \leq H_{\text{games}}$ and $\sup_{x \in \mathcal{X}} \|x \|_{\infty} \leq H_{\text{games}}$. In addition, $\sup_{\theta \in \Theta} \| \theta \|_1 \leq H_{\text{games}}$.
 
    \item \label{index-games:nuisance} There exists an estimator $\hat{p}(x)$ obeying  Assumption \ref{ass:nuisance}  with $\pi_n$ rates in either $\ell_{\infty}$ or $\ell_2$ norm.

\end{enumerate}

\end{assumption}

\begin{remark}[First-Stage Estimator of $h_0(x)$]
\label{rm:fsvee:games}
For each $i \in J_{k+1}$,  define the estimator $\widehat{h}(X_i)=\widehat{h}_k(X_i)$ 
\begin{align}
\label{eq:vee:games}
\widehat{h}(X_i)=\widetilde{\Delta}_k \cdot G_k' ( X_i' \check{\alpha}_k + \check{\Delta}_k \cdot \widehat{p}(X_i) ). 
\end{align}
where the preliminary estimator $\check{\theta} =(\check{\alpha}, \check{\Delta})$ is based on the loss $\ell_{\text{pre}}(w, t, \gamma_1)$ as in \eqref{eq:losshardgames:pre}.

\end{remark}

\begin{corollary}[Games of Incomplete Information]
\label{cor:games}
Suppose  Assumption  \ref{ass:selection} holds. Let the preliminary loss $\ell_{\text{pre}}(w, t, \gamma_1)$ be as in \eqref{eq:losshardgames:pre}  and the first-stage parameter estimate $\widehat{p}(x)$. Let the final loss be as in  \eqref{eq:losshardgames} and $\widehat{g} (x)= (\widehat{p}(x), \widehat{h}(x))$, where $\widehat{h}(x)$ as defined in  \eqref{eq:vee:games}.  Then, the statement of Theorem \ref{cor:non-ortho-rate} holds for each step of  Definition \ref{def:twostep}, where the first-step nuisance rate is $O(\pi_{n})$ and the second-step nuisance rate  $g_n =  O \left(k \left( \pi_n +\sqrt{\frac{\log p}{n}} \right) \right)$.

\end{corollary}

\section{Empirical Application}
\label{sec:emp}

In this section, we study the short-term impact of Connecticut's Jobs First welfare reform experiment on women's labor supply and welfare participation decisions.  Jobs First, a welfare-to-work assistance program, was introduced in Connecticut in the late 1990s as an alternative to the federal  Aid to Families with Dependent Children (AFDC) program. Imposing  revealed preference restrictions, \cite{KT}'s nonparametric bounds
show that Jobs First induced many women to work but led  some others to reduce their earnings in
order to receive assistance.  To gain more insight into the question, we postulate a partially linear logistic model for women's welfare participation decision and estimate  heterogeneous Jobs First  effects.

The data for our analysis are the same as in \cite{KT}. They   come from Manpower Demonstration Research Corporation (MDRC). In 1996, MDRC conducted a randomized trial that randomly selected a set of eligible female applicants into Jobs First (the treatment group), leaving the remaining females eligible for AFDC (the control group). The  outcome of interest is the binary indicator that is equal to one if a woman receives any type of welfare (e.g., AFDC or food stamps) at the fourth quarter after random assignment. Baseline characteristics  include constant, age, education level,  quarterly history of employment, earnings, AFDC, and food stamps  during $8$ quarters before random assignment (RA). We postulate a logistic specification with a partially linear index
\begin{align}
\label{eq:empspec}
\Pr (Y=1 | D=d, X=x) &= G (d \cdot ((1,x)'\theta_0) + f_0(x)),
\end{align}
where $D=1$ if a woman is assigned to Jobs First and $Y=1$ if a woman is on welfare.  We assume that the treatment $D$ and the controls $X$ affect the outcome $Y$ via an additively separable index
\begin{equation} \label{eq:index} D \cdot ((1,X)'\theta_0) + f_0(X) \end{equation}
entering the logistic link function $G(\cdot): \R \rightarrow \R$.  The base treatment $D$ enters the index  \eqref{eq:index} through its interactions with the  control vector $(1,X)$.   In addition, the control vector $X$ enters the index \eqref{eq:index} through an unknown function $f_0(x)$, which summarizes the confounding effect of the controls on $Y$.  Our main object of interest is the heterogeneous treatment effects vector 
 $$\theta_0= (\theta_{0,1}, \theta_{0,-1})',$$
where $\theta_{0,1}$ is the baseline treatment effect and $\theta_{0,-1}$ is the treatment interaction effect. We assume that, out of $20$ treatment effect interactions,  only few have non-zero value, but  do not know their identities.

A standard approach to this problem is to require the confounding function $f_0(x)$ to be a  linear sparse function of the controls
\begin{align}
\label{eq:linearity}
 f_0(x) = B(x)'\alpha_0, \quad p_{\alpha} = \dim (\alpha_0)  \gg n, \quad \| \alpha_0 \|_0 = k_{\alpha} \ll n,
\end{align}
where $B(x)$ is a vector of basis functions of $x$ and $n$ is the sample size. We take $B(x)$ to be a vector of  $p_{\alpha} = 1, 600$ pairwise interactions of the controls.  The  joint estimate of $\theta_0$ and $\alpha_0$  is 
\begin{align}
\label{eq:naive}
&(\check{\theta}_{\text{direct}}, \widetilde{\alpha}_{\text{direct}}) = \\
&\arg \min_{(\theta, \alpha) \in \R^{p + p_{\alpha}} } \dfrac{1}{n} \sum_{i=1}^n -Y_i (D_i \cdot ((1,X_i)'\theta) + B(X_i)'\alpha) + \log (1+ \exp( D_i \cdot ((1,X_i)'\theta) + B(X_i)'\alpha )) \\ \nonumber
&+\lambda_{\text{direct}} \| \theta_{-1} \|_1 + \lambda_{\text{direct}} \| \alpha \|_1, \nonumber
\end{align}
where  the penalty parameter $\lambda_{\text{direct}}$ recommended by \cite{belloni2016postselection} is
$$
\lambda_{\text{direct}} = \dfrac{1.1}{2\sqrt{n}} \Phi^{-1}(1-0.05/  (p+p_{\alpha}) \log (n) ) =  0.036.
$$
Following \cite{tan2018modelassisted}, we do not penalize the first coordinate $\theta$, imposing the sparsity assumption on the heterogeneous modification effects but not the level of the treatment effect for the baseline category.

We compare the direct estimator to the Regularized $M$-Estimator $\widehat{\theta}$, described in Section \ref{sec:nte1}. This estimator no longer requires the control function $f_0(x)$ to be sparse with respect to the chosen basis $B(x)$. Instead, we assume that  the conditional probability of welfare 
$$
G_0(d,x)= \E[ Y | D=d, X=x] 
$$
can be well-approximated by trees.  While this assumption is less interpretable, it allows $G_0(d,x)$ to have a nonlinear argument. We estimate  $G_0(d,x)$ by the probability random forest as in \cite{Malley} with $100$ trees and default size of leaf node. As for $q_0(x)$, we take
$$
\widehat{q}(x) = \widehat{\E} \bigg[ \log \dfrac{\widehat{G}(d,x)}{1-\widehat{G}(d,x)} \bigg| X=x \bigg],
$$
where the outer expectation function is estimated by regular random forest of \cite{Breiman}. In addition,  we assume that the propensity score $p_0(x)$  in equation \eqref{eq:cexp}  is a sufficiently  smooth function of $x$ and estimate it by simple logistic regression.  The weighting function $\widehat{V}(d,x)$ is estimated as in \eqref{eq:vee}. In Definition \ref{def:alg:lasso}, we set $\lambda$ to be
$$
\lambda=\lambda_{\text{ortho}} = \dfrac{1.1}{2\sqrt{n}} \Phi^{-1}(1-0.05/n ) =  0.031,
$$
standardizing each covariate after interacting it with treatment. The Lasso estimate $\widehat{\theta}$  and its post-penalized analog post-Lasso-logistic $\widehat{\theta}_{\text{PL}}$ are reported in Table \ref{tab:res}, Columns (2)-(3). Finally,  we report the unpenalized version of the orthogonal estimator defined as
\begin{align}
\label{eq:nopenalty}
\widehat{\theta}_{\text{unpenalized}}= \arg \min_{\theta \in \mathrm{R}^p} \widehat{Q}(\theta, \widehat{g}).
\end{align}
Since the number of treatment interactions $p=20$ is less than the sample size, this estimator is well-defined. We report this estimate and the standard errors in Table \ref{tab:res}, Columns (4)-(5). Our method is straightforward to implement using the \url{glmnet} function in the  \url{glmnet} $R$ package, setting the \url{offset} argument to $\widehat{q}(x)$, the \url{weights} argument to $\widehat{V}(d,x)$, and the 
\url{penalty.factor}=$(0,1,1,\dots, 1)$ to accommodate  the penalty in \eqref{eq:naive}\begin{footnote}{The package is available at \url{https://github.com/vsyrgkanis/plugin_regularized_estimation}}\end{footnote}.

 Our empirical findings are as follows. The direct estimator (Column (1)) implies that Jobs First unambiguously increased the fraction of women on welfare for all women. In contrast, the orthogonal estimator (Column (2)) implies that women with long history of prior AFDC receipt (at least 5 months) were pushed out of the assistance. This pattern makes sense. The Jobs First (treatment group) faced a time limit of 21 months on welfare while the AFDC (control group) had no time limit. The orthogonal estimator picks up this pattern, while the direct one fails to do so. Finally, the unpenalized model in Columns (4) and (5) finds evidence in favor of treatment effect heterogeneity, but most interaction effects are not significant at  $\alpha=0.05$. Therefore, orthogonal Lasso captures a more nuanced pattern of treatment effect that is left out by direct Lasso due to a substantially larger complexity of the nuisance parameter relative to the target one,  failure of the sparsity assumption \eqref{eq:linearity}, or both.

\begin{table}
\caption{ Jobs First heterogenous effects on welfare participation}
\begin{center}
\begin{tabular}{ cccccc}
 \toprule
 & (1) & (2) & (3) & (4) & (5) \\
 & direct Lasso & ortho Lasso & post-ortho-Lasso  & not penalized & not penalized \\
 & $\check{\theta}_{\text{direct}}$ & $\widehat{\theta}$ & $\widehat{\theta}_\text{PL}$ & $\widehat{\theta}_{\text{unpenalized}}$ & SE($\widehat{\theta}_{\text{unpenalized}}$) \\
\toprule
intercept & 0.301 & 0.528 & 0.976 & 0.035 & 2.056 \\ 
  white &  &  &  & 0.906 & 1.427 \\ 
  black &  &  &  & 1.023 & 1.456 \\ 
  hisp &  &  &  & 0.342 & 1.472 \\ 
  marnvr &  &  &  & -0.395 & 1.591 \\ 
  marapt &  &  &  & -0.167 & 1.607 \\ 
  yrern &  &  &  &  0.000&0.000  \\ 
  yrernsq &  &  &  &0.000 &  0.000\\ 
  yradc &  &  &  &  0.000&0.000  \\ 
  yrfst &  &  &  &   0.001 &0.000\\ 
  yremp &  &  &  & -0.003 & 0.011 \\ 
  yrvad &  &  &  & 0.006 & 0.022 \\ 
  yrvfs &  &  &  & -0.002 & 0.020 \\ 
  yrkvad &  & -0.022 & -0.136 & -0.115 & 0.172 \\ 
  anyernpq &  &  &  & 0.061 & 1.087 \\ 
  anyadcpq &  &  &  & 0.212 & 1.940 \\ 
  anyfstpq &  &  &  & -0.511 & 1.938 \\ 
  nohsged &  &  &  & 0.575 & 0.612 \\ 
  applcant &  &  &  & 0.182 & 0.899 \\
 
 \bottomrule
\end{tabular}
\end{center}
\caption*{Notes. The data come from \cite{KT}. Table reports heterogeneous Jobs First effects in the equation  \eqref{eq:empspec}, estimated by the direct Lasso, ortho Lasso, post-ortho-Lasso-logistic, and unpenalized approaches. Each row corresponds to a treatment interaction with a covariate. The selected covariate  \text{yrkvad} shows the number of months on AFDC in the year prior RA. The other covariate names are: $\text{white}, \text{black}, \text{hisp}$ are race categories; marnvr is never married, marapt is widowed or separated; $\text{yrern}/\text{yrernsq}/\text{yradc}/\text{yrfst}$  the total earnings/earnings squared/AFDC/food stamps in dollars  the year prior to RA; $\text{yremp}/\text{yrvad}/\text{yrvfs}$ is a binary indicator equal to one   if a woman received  any earnings/AFDC/food stamps in the year before RA; $\text{anyernpq}/\text{anyadcpq}/\text{anyfstpq}$ are binary indicators equal to one   if a woman received  any earnings/AFDC/food stamps in her lifetime; $\text{nohsged}=1$ if the subject has no high school or GED degree and $0$ otherwise; $\text{applcant}=1$ if female applied to Jobs First and $0$ otherwise. Estimators: $\check{\theta}_{\text{direct}}$ as in \eqref{eq:naive} under assumption \ref{eq:linearity}, the Lasso estimator $\widehat{\theta}$ of Definition \ref{def:twostep}, post-Lasso-logistic  $\widehat{\theta}_\text{PL}$, and unpenalized logistic estimator $\widehat{\theta}_{\text{unpenalized}}$ defined in equation \eqref{eq:nopenalty}.
The sample size $n=4, 642$ females.  For the estimation details, see text. 
  }
\label{tab:res}
\end{table}
\newpage

\renewcommand{\theequation}{S.\arabic{equation}}
\renewcommand{\thesection}{S\arabic{section}}
\renewcommand{\thetheorem}{S.\arabic{theorem}}
\renewcommand{\theassumption}{S.\arabic{assumption}}
\renewcommand{\theproposition}{S.\arabic{proposition}}
\renewcommand{\thecorollary}{S.\arabic{corollary}}
\renewcommand{\thelemma}{S.\arabic{lemma}}
\renewcommand{\theexample}{S.\arabic{example}}
\renewcommand{\theremark}{S.\arabic{remark}}

\setcounter{equation}{0}
\setcounter{theorem}{0}
\setcounter{lemma}{0}
\setcounter{assumption}{0}

\renewcommand{\thesection}{A}

\section{General Theory of Extremum Estimators}
\label{sec:theory}

In this section, we introduce Regularized Extremum Estimator and establish its properties.  Suppose there exists a population loss  $Q(\theta, g)$  and its sample analog $\widehat{Q}(\theta, \widehat{g})$. Define the target parameter $\theta_0$  as the unique minimizer of $Q(\theta, g_0$)  (i.e.,  \eqref{eq:ld} holds), and the Regularized Extremum Estimator as in \eqref{eq:alg:lasso}. In contrast to \eqref{eq:sampleloss},  $\widehat{Q}(\theta, \widehat{g})$ may not be a sample average. Finally, let the set $\mathbb{B}$ be as in  \eqref{eq:mathbb} and the event $\mathcal{V}$ be as in \eqref{eq:event}.

\begin{assumption}[Convexity]\label{ass:convloss}
With probability one, for any $g$ in the set $\mcG_n$, the function $\theta \rightarrow \widehat{Q}(\theta,g)$ is a convex twice differentiable function of $\theta$ on an open convex set that contains the parameter space $\Theta$.
\end{assumption}

\begin{assumption}[Restricted Strong Convexity on $\mathbb{B}$]
\label{ass:boundedhess}
There exists a constant $\bar{\gamma}>0$ so that the curvature of $Q(\theta, g_0)$ defined in \eqref{eq:rsc} is bounded from below by $\bar{\gamma}$:
\begin{align*}
\inf_{\theta \in \mathbb{B}, \nu=\theta - \theta_0 }  \dfrac{\nu^T  \nabla_{\theta \theta} Q(\theta, g_0)  \nu} { \| \nu \|_2^2}  \geq \bar{\gamma}.
\end{align*}
\end{assumption}

\begin{assumption}[Uniform Convergence on $\mathbb{B}$]
\label{ass:grc}
There exists a sequence $\tau_n=o(1)$ so that $\nabla_{\theta \theta} \widehat{Q}(\theta, g_0) $ uniformly converges  to $\nabla_{\theta \theta} Q(\theta,g_0)$:
\begin{align}
\sup_{\theta \in \mathbb{B}} \| \nabla_{\theta \theta} \widehat{Q}(\theta, g_0) - \nabla_{\theta \theta} Q(\theta, g_0) \|_{\infty} = O_{P} (\tau_n).
\end{align}

\end{assumption}

\begin{assumption}[Uniformly Lipschitz Hessian on $\mathbb{B}$]
\label{ass:liphessian-non-conv}
For any $g$ and $g'$ in $\mcG_n$,  there exists a sequence $\xi_n=o(1)$ so that  the following bound holds
\begin{align}
\sup_{\theta \in \mathbb{B}} \| \nabla_{\theta \theta} \widehat{Q}(\theta, g) - \nabla_{\theta \theta} \widehat{Q}(\theta, g') \|_{\infty} = O_{P}(g_n+ \xi_n).
\end{align}
\end{assumption}

\begin{assumption}[Convergence Rate of Empirical Gradient]\label{ass:lossconv}
For any  $g\in \mcG_n$, there exists a sequence $\epsilon_n$ such that $ \|\nabla_{\theta} \widehat{Q}(\theta_0, g) - \nabla_{\theta} Q(\theta_0, g)\|_{\infty}=O_{P}(\epsilon_n)$.
\end{assumption}

We describe the influence of the estimation error $\widehat{g}-g_0$ on the gradient  $\nabla_{\theta} Q(\theta_0, \cdot)$ using  pathwise derivatives w.r.t the nuisance parameter $g_0$. The first-order derivative is
\begin{equation}
    D_{r}[g-g_0, \nabla_{\theta}Q(\theta_0, g )]:=\dfrac{\partial}{\partial r}\nabla_{\theta}Q(\theta_0, r(g-g_0)+g_0)
\end{equation}
and the second-order derivative is
\begin{equation}
    D_{r}^2[g-g_0, \nabla_{\theta}Q(\theta_0, g )]:=\dfrac{\partial^2}{\partial r^2}\nabla_{\theta}Q(\theta_0, r(g-g_0)+g_0).
\end{equation}

\begin{assumption}[Bounded  Gradient of Population Loss w.r.t. Nuisance]\label{ass:nonortholoss}
For any $g \in \mcG_n$, there exist constants $B_0$ and $B$ so that
$$ \forall r \in [0,1):  \left\| D_0 [g-g_0, \nabla_{\theta}Q(\theta_0, g_0 )] +  D_{r}^2[g-g_0, \nabla_{\theta}Q(\theta_0, g_0 )] \right\|_{\infty} \leq B_0 \| g-g_0\| + B \|g - g_0\|^2.  $$
\end{assumption}

\begin{theorem}[Regularized Extremum Estimator]\label{thm:main-convex}  Suppose Assumptions \ref{ass:nuisance} and \ref{ass:convloss}-\ref{ass:nonortholoss}  hold with $k(\tau_{n} + g_{n} + \xi_{n})=o(1)$.  For $\lambda = 2 C(\epsilon_{n} + B_0 g_n + B g_n^2)$, and for $C$ and $n$ large enough,

\begin{align}
    \| \widehat{\theta} - \theta_0 \|_2 &\lesssim_P \sqrt{k}  (\epsilon_n  + B_0 g_n +  g_n^2) &
     \| \widehat{\theta} - \theta_0 \|_1 &\lesssim_P  k (\epsilon_n  + B_0 g_n +  g_n^2) .
\end{align}

\end{theorem}

Theorem \ref{thm:main-convex} establishes the convergence rate for Regularized  Extremum Estimator. Its building blocks are described in the following lemmas.

\begin{lemma}[Convexity and Restricted Subspace, \cite{Negahban}] If $\frac{\lambda}{2} \geq \|\nabla_\theta \widehat{Q}(\theta_0, \widehat{g})\|_{\infty}$ and $\theta \rightarrow Q(\theta, \widehat{g})$ is a convex function of $\theta$, then $\nu \in {\cal C}(T; 3)$, where $\nu=\widehat{\theta}-\theta_0$.\label{lem:subspace}
\end{lemma}

\begin{proof}[Proof of Lemma \ref{lem:subspace}]

By definition of $\widehat{\theta}$,
\begin{equation}\label{eqn:minimizer}
\widehat{Q}(\widehat{\theta},\widehat{g}) - \widehat{Q}(\theta_0,\widehat{g}) \leq \lambda \left(\|\theta_0\|_1 - \|\widehat{\theta}\|_1\right).
\end{equation}
Additivity of $\ell_1$ norm and triangle inequality imply
\begin{align*}
\|\widehat{\theta}\|_1 &= \|\theta_0 + \nu_T\|_1 + \|\nu_{T^c}\|_1 \geq \|\theta_0\|_1 - \|\nu_T\|_1 + \|\nu_{T^c}\|_1  \\
\lambda \left(\|\theta_0\|_1 - \|\widehat{\theta}\|_1\right) &\leq \lambda \left(\|\nu_T\|_1 - \|\nu_{T^c}\|_1\right).
\end{align*}
Convexity of $\theta \rightarrow \widehat{Q}(\theta,\widehat{g})$, Cauchy-Schwarz inequality, and the choice of $\lambda$ imply:
\begin{equation}\label{eqn:convexity}     
\widehat{Q}(\widehat{\theta}, \widehat{g}) - \widehat{Q}(\theta_0, \widehat{g}) \geq \nabla_\theta \widehat{Q}(\theta_0, \widehat{g}) \cdot (\widehat{\theta}-\theta_0) \geq -\|\nabla_\theta \widehat{Q}(\theta_0, \widehat{g})\|_\infty \|\nu\|_1 \geq - \frac{\lambda}{2} \|\nu\|_1.
\end{equation}
Combining equations \eqref{eqn:minimizer} and \eqref{eqn:convexity} gives:
\begin{equation*}
\lambda \left(\|\nu_T\|_1 - \|\nu_{T^c}\|_1\right) \geq \widehat{Q}(\widehat{\theta},\widehat{g}) - \widehat{Q}(\theta_0,\widehat{g}) \geq -\frac{\lambda}{2} \|\nu\|_1.
\end{equation*}
Dividing by $\lambda$ and re-arranging gives $3\|\nu_T\|_1 \geq \|\nu_{T^c}\|_1$.
\end{proof}

\begin{lemma} [Restricted Strong Convexity for Empirical Loss]  
\label{lem:fsrsc}
Suppose  Assumptions \ref{ass:convloss} and \ref{ass:boundedhess} hold.   On the event $\mathcal{V}$, the following bound holds:
\begin{align*}
 \inf_{  \nu \in {\cal C}(T;3)}     \frac{ (\nabla_\theta  \widehat{Q}(\theta_0 + \nu, \widehat{g})  - \nabla_\theta  \widehat{Q}(\theta_0, \widehat{g})) \cdot \nu }{ \| \nu \|_2^2}  \geq \bar{\gamma}/2.
\end{align*}

\end{lemma}
\begin{proof}[Proof of Lemma \ref{lem:fsrsc}]
\textbf{ Step 1. } For some $\bar{r} \in [0,1)$ that may depend on $\widehat{g}$, mean value theorem implies
\begin{align*}
	   ( \nabla_\theta  \widehat{Q}(\theta_0 + \nu, \widehat{g})  - \nabla_\theta  \widehat{Q}(\theta_0, \widehat{g})) \cdot \nu &= \nu^T \cdot  \nabla_{\theta \theta}  \widehat{Q}(\theta_0 + \bar{r} \nu, \widehat{g}) \cdot \nu. 
\end{align*}
Invoking Cauchy-Schwartz, definition of  ${\cal C}(T; 3)$ in \eqref{eq:cone}, and assumption of the theorem gives:
\begin{align*}
&\sup_{\bar{r} \in [0,1), \nu \in {\cal C}(T;3)} \frac{| \nu^T \cdot  ( \nabla_{\theta \theta} \widehat{Q}(\theta_0 + \bar{r} \nu, \widehat{g})   -   \nabla_{\theta \theta}  Q(\theta_0 + \bar{r} \nu, g_0) )   \cdot \nu  |}{\|\nu\|^2_2} \\
&\leq   \sup_{\nu \in {\mathcal C}(T;3)}\frac{\|\nu \|_1^2 \bar{\gamma}/(32k) }{\|\nu\|^2_2}  \tag{Cauchy-Schwartz} \\
 &\leq \sup_{\nu \in {\mathcal C}(T;3)} \frac{((1+3) \| \nu_T \|_1)^2 \bar{\gamma}/(32k)   }{\|\nu\|^2_2}  \tag{Definition of ${\cal C}(T; 3)$} \\
 &\leq 16 k \bar{\gamma}/(32k)   < \bar{\gamma}/2.
\end{align*}

\textbf{ Step 2. } Triangle inequality implies
\begin{align*}
   & \inf_{  \nu \in {\cal C}(T;3)}     \frac{   \nu^T \cdot  \nabla_{\theta \theta}  \widehat{Q}(\theta_0 + \bar{r} \nu, \widehat{g}) \cdot \nu  }{ \| \nu \|_2^2} \\
   &\geq  \inf_{  \nu \in {\cal C}(T;3)}     \frac{   \nu^T \cdot  \nabla_{\theta \theta}  Q(\theta_0 + \bar{r} \nu, g_0) \cdot \nu  }{ \| \nu \|_2^2} - 
   \sup_{  \nu \in {\cal C}(T;3)}     \frac{ |  \nu^T \cdot ( \nabla_{\theta \theta} \widehat{Q}(\theta_0 + \bar{r} \nu, \widehat{g})   -   \nabla_{\theta \theta}  Q(\theta_0 + \bar{r} \nu, g_0) ) \cdot \nu  |}{ \| \nu \|_2^2} \\
   &\geq \bar{\gamma} - \bar{\gamma}/2 =\bar{\gamma}/2.
\end{align*}

\end{proof}

 \begin{lemma}[Oracle Inequality]\label{thm:gen-conv}
Suppose  Assumptions \ref{ass:convloss} and \ref{ass:boundedhess} hold.     On the intersections of events \eqref{eq:eventnoise} and   \eqref{eq:event},  the following inequality holds:
\begin{align}
\|\widehat{\theta} - \theta_0\|_2 \leq& \frac{3\sqrt{k}}{\bar{\gamma} } \lambda  & \|\widehat{\theta} - \theta_0\|_1 \leq&  \frac{12 k}{\bar{\gamma} } \lambda. 
\end{align}
\end{lemma}

\begin{proof}[Proof of Lemma \ref{thm:gen-conv}]
 Let $\nu = \widehat{\theta} - \theta_0$.  Invoking Lemmas  \ref{lem:subspace}  and \ref{lem:fsrsc} gives:
 \begin{align*}
 \lambda \left(\| \theta_0 \|_1 - \|\widehat{\theta}\|_1 \right)   \geq~& \nabla_\theta \widehat{Q}(\widehat{\theta}, \widehat{g}) \cdot (\widehat{\theta}-\theta_0) \tag{Optimality of $\widehat{\theta}$} \\
=~&  \nabla_\theta \widehat{Q}(\theta_0, \widehat{g}) \cdot (\widehat{\theta}-\theta_0) +  ( \nabla_\theta  \widehat{Q}(\widehat{\theta}, \widehat{g})  - \nabla_\theta  \widehat{Q}(\theta_0, \widehat{g})) \cdot (\widehat{\theta}-\theta_0)   \\
\geq~& \nabla_\theta \widehat{Q}(\theta_0, \widehat{g}) \cdot (\widehat{\theta}-\theta_0)  + \frac{\bar{\gamma}}{2} \|\nu\|_2^2  \tag{Lemmas \ref{lem:fsrsc}, \ref{lem:subspace}}\\
\geq~& -\|\nabla_{\theta} \widehat{Q}(\theta_0,\widehat{g})\|_{\infty} \cdot \|\nu\|_1  + \frac{\bar{\gamma}}{2} \|\nu\|_2^2  \tag{Cauchy-Schwarz}\\
\geq~& -\frac{\lambda}{2} \cdot \|\nu\|_1  + \frac{\bar{\gamma}}{2} \|\nu\|_2^2 \tag{Assumption on $\lambda$}
\end{align*}
Rearranging the inequality gives
\begin{align*}
\frac{\bar{\gamma}}{2} \|\nu\|_2^2 \leq \frac{3\lambda}{2} \|\nu_T\|_1 - \frac{\lambda}{2} \|\nu_{T^c}\|_1 \leq \frac{3\lambda}{2} \|\nu_T\|_1 \leq \frac{3\lambda \sqrt{k}}{2} \|\nu\|_2.
\end{align*}
Dividing over by $\|\nu\|_2$ yields the theorem. By Lemma \ref{lem:subspace}, $\nu \in {\cal C}(T;3)$ and $\|\nu\|_1\leq 4 \sqrt{k} \|\nu\|_2$.
\end{proof}

\begin{lemma}[Influence on Oracle Gradient]\label{lem:influence}
Suppose Assumptions \ref{ass:nuisance},  \ref{ass:lossconv} and \ref{ass:nonortholoss} hold. With probability $1-o(1)$,
\label{lem:oraclegrad}
\begin{align}
    \label{eq:grad}
     \|\nabla_\theta \widehat{Q}(\theta_0, \widehat{g})\|_{\infty} &=  O_{P} (\epsilon_{n} + B_0 g_{n} + B g_{n}^2).
\end{align}
\end{lemma}

\begin{proof}[Proof of Lemma \ref{lem:influence}]
In what follows, we condition on the event $\mathcal{E}_n$, which holds with probability $1-o(1)$. Triangle inequality implies:
    \begin{align}
        \|\nabla \widehat{Q}(\theta_0, \widehat{g})\|_{\infty} &\leq  \|  \nabla_{\theta} [Q(\theta_0, \widehat{g}) -Q(\theta_0, g_0) ] \|_{\infty} + \|  \nabla_{\theta} [\widehat{Q}(\theta_0, \widehat{g}) -Q(\theta_0, \widehat{g})] \|_{\infty} \label{eqn:gradient-eqn-1} \\
        &= I + II \nonumber.
    \end{align}
     Conditional on the event $\mathcal{E}_n$, Assumption \ref{ass:lossconv} implies
      $$II=\|  \nabla_{\theta} [\widehat{Q}(\theta_0, \widehat{g}) -Q(\theta_0, \widehat{g})] \|_{\infty} = O_{P} (\epsilon_n).$$
     By Lemma 6.1 from \cite{chernozhukov2016double}, this statement holds unconditionally.  
   On the event $\mathcal{E}_n$, by Assumption \ref{ass:nonortholoss}
    \begin{align*}
  I= \|  \nabla_{\theta} [Q(\theta_0, \widehat{g}) - Q(\theta_0, g_0) ] \|_{\infty}   &\leq  \sup_{g \in \mcG_n} \|  \nabla_{\theta} [Q(\theta_0, g) - Q(\theta_0, g_0) ] \|_{\infty} \\
   &\leq B_0 g_{n}+ B g_{n}^2. 
    \end{align*}
     Thus, \eqref{eq:grad} follows.

\end{proof}

\begin{proof}[Proof of Theorem \ref{thm:main-convex}]
 By Lemma \ref{lem:oraclegrad} and assumptions of the Theorem, 
\begin{align}
\label{eq:firsthalf}
  \lambda/2 \geq  \| \nabla_\theta \widehat{Q}(\theta_0, \widehat{g}) \|_{\infty}  \text{ holds  w.p. } 1-o(1).
\end{align} 
Decomposing the difference $\nabla_{\theta \theta} \widehat{Q}(\theta, \widehat{g} )  - \nabla_{\theta \theta} Q(\theta, g_0 )$ gives
\begin{align*}
\sup_{\theta \in \mathbb{B}} \| \nabla_{\theta \theta} \widehat{Q}(\theta, \widehat{g} )  - \nabla_{\theta \theta} Q(\theta, g_0 )  \|_{\infty}  &\leq \sup_{\theta \in \mathbb{B}}  \| \nabla_{\theta \theta} \widehat{Q}(\theta, \widehat{g} )  - \nabla_{\theta \theta} \widehat{Q}(\theta, g_0 ) \|_{\infty} \\
&+ \sup_{\theta \in \mathbb{B}}  \| \nabla_{\theta \theta} \widehat{Q}(\theta, g_0)  - \nabla_{\theta \theta} Q(\theta, g_0 ) \|_{\infty}\\
&= A+ B.
\end{align*}
By Assumptions \ref{ass:nuisance}, \ref{ass:liphessian-non-conv}, \ref{ass:grc}, on the event $\mathcal{E}_n$,
\begin{align*}
A &=  O_P(g_n + \xi_n), \quad B = O_P(\tau_n).
\end{align*}
By assumption of the Theorem, $k (g_n + \xi_n + \tau_n) = o(1)$. Therefore, for $n$ large enough, 
\begin{align}
\sup_{\theta \in \mathbb{B}}  \| \nabla_{\theta \theta} \widehat{Q}(\theta, \widehat{g} )  - \nabla_{\theta \theta} Q(\theta, g_0 )  \|_{\infty}  \leq \bar{\gamma}/(32k)  \text{ holds  w.p. } 1-o(1).
\end{align}

\end{proof}

\renewcommand{\thesection}{B}

\section{Proofs of Section 3}

\label{sec:proofs_moment}

In this section, we prove Theorem  \ref{cor:non-ortho-rate} as a special case of Theorem \ref{thm:main-convex}. We also prove   Theorem \ref{lem:general}.   Define an event $\mathcal{E}_n := \{ \widehat{g}_{k} \in \mcG_n \quad \forall k \in [K] \}$, such that the nuisance parameter estimate $\widehat{g}_k$ belongs to the realization set $\mcG_n$ for each fold $k \in [K]$. By union bound, this event holds w.h.p. $$\Pr (\mathcal{E}_n ) \geq 1- K \epsilon_n = 1-o(1).$$ 
For a given partition $k$ in $\{1,2, \dots, K\}$, define the partition-specific averages $$\Enk f(W_i) := \dfrac{1}{n_k} \sum_{i \in J_k} f(W_i), \quad \Gnk f(W_i) := \dfrac{1}{\sqrt{n_k}} \sum_{i \in J_k} [f(W_i) - \int f(w) dP (w)] .$$

\begin{lemma}[Verification of Assumption \ref{ass:convloss}]
\label{lem:obvious}
Assumption \ref{ass:convexitymoment}  implies Assumption \ref{ass:convloss}.
\end{lemma}
\begin{proof}[Proof of Lemma \ref{lem:obvious}]
If $m(w,t,\gamma)$ is non-decreasing in $t$,  the loss $\ell(w,t,\gamma)$ is convex in $t$. Plugging the linear function
 $\theta \rightarrow t=\Lambda(z, g)'\theta$ into the  loss sketch $\ell(w,t,\gamma)$ gives the convex sample loss $\widehat{Q} (\theta, g)$.

\end{proof}

\begin{lemma}[Verification of Assumption \ref{ass:boundedhess}]
\label{lem:boundedhess}
Assumptions \ref{ass:identification} with $C_{\text{min}}$  and  \ref{ass:identification2} with $B_{\text{min}}$ imply Assumption \ref{ass:boundedhess} with $\bar{\gamma} = B_{\text{min}} C_{\text{min}}$. 
\end{lemma}

\begin{proof}[Proof of Lemma \ref{lem:boundedhess}]
Observe that
\begin{align*}
&\inf_{\theta \in \mathbb{B} }  \nu^\top  \nabla_{\theta \theta} Q(\theta, g_0) \nu  \\
=&\inf_{\theta \in \mathbb{B} }  \E \bigg[ \nabla_t m(W, t, \gamma) \bigg|_{\gamma = g_0(Z), t=\Lambda(Z, \gamma)'\theta} | Z=z \bigg] \min \eig   \, \Sigma \, \| \nu \|_2^2 \\
&\geq B_{\text{min}} C_{\text{min}}  \| \nu \|_2^2.
\end{align*}

\end{proof}

\begin{lemma}[Verification of Assumption \ref{ass:grc}]
Lemma \ref{lem:grc} verifies Assumption \ref{ass:grc}.

\end{lemma}

\begin{lemma}[Verification of Assumption \ref{ass:liphessian-non-conv}]
\label{lem:liphessian-non-conv}
Assumptions \ref{ass:nuisance} and \ref{ass:smooth}  with a constant $U$ implies   Assumption \ref{ass:liphessian-non-conv} with $\xi_n = 3  d \cdot U^3 \, \sqrt{\frac{\log(2p)}{n}}$ for $\mathbb{B}$ as in \eqref{eq:mathbb} and $\| g \|_r$ as in  \eqref{eq:norms} for $r \in \{\infty,2\}$.

\end{lemma}

\begin{proof}[Proof of Lemma \ref{lem:liphessian-non-conv}]
\textbf{ Step 1. } Observe that  for each $(v,j) \in \{1,2,\dots, p\}^2$,
\begin{align}
\nabla_{ \theta_v \theta_j} \ell(w, \Lambda(z, \gamma)'\theta, \gamma )=\nabla_{ t} m(w, \Lambda(z,\gamma)'\theta, \gamma) \cdot \Lambda_v(z, \gamma) \cdot \Lambda_j(z, \gamma) 
\end{align}
has its gradient in $\gamma$ bounded by  $L=3 \cdot U^3 $ in the absolute norm.

Verification of Assumption \ref{ass:liphessian-non-conv} for $\| g \| = \| g \|_{\infty}$. Assumptions \ref{ass:smooth} implies Assumption \ref{ass:liphessian-non-conv} with $L=U^3$.
\begin{align*}
   &|\nabla_{\theta_v\theta_j} \widehat{Q}(\theta, g) - \nabla_{\theta_v\theta_j} \widehat{Q}(\theta, g')|   \leq 3 \cdot U^3 \sup_{z\in \mcZ}   \|g(z)-g'(z)\|_1.   \end{align*}

\textbf{ Step 2. } Verification of Assumption \ref{ass:liphessian-non-conv} for $\| g \| = \| g \|_2$. For each $(v,j) \in \{1,2,\dots, p\}^2$, the function  above has $U$-bounded derivative, its derivative is bounded by $L=3U^3$.
By McDiarmid's inequality, for any fixed $g, g'\in \mcG_n$, with probability $1-o(1)$,
\begin{align*}
    \left|[\Enk-\E][\|g(Z)-g'(Z)\|_1] \right| \leq  3  d \cdot U^3 \, \sqrt{\frac{\log(2p)}{n_k}} = O\left( 3  d \cdot U^3 \, \sqrt{\frac{\log(2p)}{n}} \right):= O(\xi_{n} ).
\end{align*}

\end{proof}

\begin{lemma}[Verification of Assumption \ref{ass:lossconv}]
\label{lem:lossconv}
Assumption \ref{ass:smooth} with a constant $U$ implies    Assumption \ref{ass:lossconv} with $\epsilon_n = U^2 \cdot \sqrt{\frac{\log(2p)}{n}}$.
\end{lemma}

\begin{proof}[Proof of Lemma \ref{lem:lossconv}]
Conditional on the event $\mathcal{E}_n$, the estimate $g=\widehat{g}$ belongs to $\mcG_n$. Observe that $$\nabla_{\theta_j} \ell(w,\Lambda(z, \gamma)' \theta_0,\gamma)|_{\gamma = g(z)} = m(w, \Lambda(z, \gamma)' \theta_0,  \gamma) \cdot  \Lambda_j(z, \gamma) |_{\gamma = g(z)}, \quad j=1,2,\dots p$$ is an a.s. bounded function of data vector $w$,  bounded in absolute value by $U^2$. By McDiarmid's inequality, the sample average of $j$'th function is within $U^2 \cdot \sqrt{\frac{\log(2p/\delta)}{n_k}}$ of its mean with probability at least $1-\frac{\delta}{p}$.  Taking a union bound over the $p$ coordinates,  we get that each coordinate is within $U^2 \cdot \sqrt{\frac{\log(2p/\delta)}{n_k}}$ from its respective mean with probability at least $1-\delta$. Therefore, $\epsilon_n$ can be taken to be $\epsilon_n = U^2 \cdot \sqrt{\frac{\log(2p)}{n_k}} = O(U^2 \cdot \sqrt{\frac{\log(2p)}{n}})$. \end{proof}

\begin{lemma}[Verification of Assumption \ref{ass:nonortholoss}]
\label{lem:orthohessian}
Assumption \ref{ass:smooth} with a constant $U$ implies Assumption \ref{ass:nonortholoss} with $B_0 = U^2$ and $B = 4U^2$.
\end{lemma}
\begin{proof}[Proof of Lemma \ref{lem:orthohessian} ]
By Assumption \ref{ass:smooth}, for any $w, \theta$ and $g \in \mcG_n$, define the  matrix $A_j \in \mathrm{R}^{d \bigtimes d}$ as
$$ A_j:= A_j(w, \theta, g) =  \nabla_{\gamma\gamma} \bigg[ m (w, t, \gamma)|_{\gamma = g(z), t = \Lambda(z, \gamma)'\theta} \cdot  \Lambda_j(z, \gamma) \bigg]. $$
By Cauchy-Schwarz inequality,
 \begin{align*}
| (g(Z)-g_0(Z))' A_j (g(Z)-g_0(Z)) | &\leq \| A_j \|_{\infty} \| (g(Z)-g_0(Z)) \|_1^2.
\end{align*} 
Taking expectations on each side and invoking  $\sup_{ 1 \leq j \leq p }\| A_j \|_{\infty}  \leq 4 \cdot U^2 \quad \text{ a.s. }$  gives
\begin{align*}
    D_r^2[g-g_0, \nabla_{\theta_j} Q(\theta_0, g_0)]    &\leq \E | (g(Z)-g_0(Z))' A_j (g(Z)-g_0(Z)) | 
     &\leq 4 \cdot U^2   \cdot  \E[\|g(Z)-g_0(Z)\|_1^2]\\
    &\leq 4 \cdot U^2  \, \sup_{z\in \mcZ} \|g(z)-g_0(z)\|_{1}^2\\
   &\leq 4 \cdot U^2  \, \|g-g_0\|_{\infty}^2,
\end{align*}
Likewise, $\| D_0[g-g_0, \nabla_{\theta} Q(\theta_0, g_0)] \|_{\infty} \leq  U^2   \E[\|g(Z)-g_0(Z)\|_1^2 $.

\end{proof}

\begin{proof}[Proof of Theorem \ref{cor:non-ortho-rate}]
Theorem \ref{cor:non-ortho-rate} follows from the statement of Theorem \ref{thm:main-convex} and Lemmas \ref{lem:obvious}-\ref{lem:orthohessian}.
\end{proof}

\begin{proof}[Proof of Theorem \ref{lem:general}]

\textbf{ Step 1. } The validity of the adjusted CMR can be seen from
\begin{align*}
\E \bigg[ m (W, \Lambda(Z, \gamma_1)'\theta_0, \gamma)|_{\gamma=g_0(Z)} \bigg| Z=z \bigg] &= \E \bigg[ m_{\text{pre}} (W, \Lambda(Z, \gamma_1)'\theta_0, \gamma_1)|_{\gamma_1=p_0(Z) } \bigg| Z=z \bigg] \\
&- h_0(z) I_0^{-1}(z) \E \bigg[ R(W, p_0(Z)) \bigg| Z=z \bigg] =0,
\end{align*}
which follows from \eqref{eq:main:pre} and \eqref{eq:cond_exog}.

\textbf{ Step 2. } We verify the orthogonality condition \eqref{eq:orthogfinal} in two steps. First, we compute the derivative of \eqref{eq:mom}  with respect to  $\gamma=(\gamma_1, \gamma_2, \gamma_3)$
\begin{align*}
&\E \bigg[ \nabla_{\gamma} m (W, \Lambda(z, \gamma_1)'\theta_0, \gamma) \cdot \Lambda(z, \gamma_1)|_{\gamma = g_0(z) }  \bigg| Z=z \bigg] \\
 + &\E \bigg[  m (W, \Lambda(z, \gamma_1)'\theta_0, \gamma) \cdot  \nabla_{\gamma}  \Lambda(z, \gamma_1) |_{\gamma = g_0(z)} \bigg| Z=z \bigg] \\
 &= (0,0,0)' + (0,0,0)' = i + ii,
\end{align*}
where $i$ is shown in Step 3 and $ii$ follows from Step 1.

\textbf{ Step 3. } The partial derivatives with respect to $\gamma=(\gamma_1, \gamma_2, \gamma_3)$ are mean zero conditionally on $z$. 
\begin{align*}
 &\E \bigg[  \nabla_{\gamma_1} m (W, \Lambda(z, \gamma_1)'\theta_0, \gamma)|_{\gamma = g_0(z) }  \bigg| Z=z \bigg]  \\
 &=    \E \bigg[  \nabla_{\gamma_1} m_{\text{pre}} (W, \Lambda(z, \gamma_1)'\theta_0, \gamma_1) \bigg| Z=z \bigg] - h_0(z) I_0(z)^{-1} I_0(z) =0.
  \end{align*} 
\begin{align*}
\E \bigg[   \nabla_{\gamma_2} m (W, \Lambda(z, \gamma_1)'\theta_0, \gamma)|_{\gamma = g_0(z) } \bigg| Z=z \bigg] =  -I_0(z)^{-1} \E[ R(W, p_0(Z)) | Z=z] =0 
\end{align*} 
\begin{align*}
 \E \bigg[  \nabla_{\gamma_3} m (W, \Lambda(z, \gamma_1)'\theta_0, \gamma)|_{\gamma = g_0(z) }  \bigg| Z=z \bigg]   =  h_0(z) I_0^{-2} (z) \E[ R(W, p_0(Z)) | Z=z] =0 
\end{align*}

\textbf{ Step 4. }  Verification of Assumptions \ref{ass:convexitymoment}-\ref{ass:smooth} for the moment function $m(w, t, \gamma)$ in \eqref{eq:rho:phi}. First, observe that
$$
 h_0(z) I_0^{-1}(z) R(W, p_0(Z))
$$
does not depend on the single index $t$. Therefore,  the moment function $m(w, t, \gamma)$ is monotone in $t$ if and only if  $m_{\text{pre}}(w, t, \gamma_1)$ is monotone in $t$. Likewise, $m(w, t, \gamma)$ is single index in $t$ if and only if  $m_{\text{pre}}(w, t, \gamma_1)$ is single index in $t$.  Finally,  $$\nabla_t m(w, t, \gamma) = \nabla_t m_{\text{pre}}(w, t, \gamma_1)$$ for any value of $w, t, \gamma$.  The function $(w, \gamma) \rightarrow \gamma_2 \gamma_3^{-1} R(w, \gamma_1)$ is a smooth function of $\gamma=(\gamma_1, \gamma_2, \gamma_3)$  by assumption of the theorem, implying Assumption \ref{ass:smooth}. 

   \end{proof}

\renewcommand{\thesection}{C}
\section{Proofs of Section 4}  
  
  Let $\check{\theta}$ be a preliminary estimator of $\theta_0$ converging at rate $\theta_n$. Lemma \ref{lem:ver:nuisance}  establishes the shrinkage rates of the nuisance realization sets introduced below.

\begin{enumerate}
\item Suppose $$ R_h(z, t,\gamma_1)= \E[\nabla_{\gamma_1} m_{\text{pre}}(W, t,\gamma_1)|Z=z]$$ is a  known function of $z, t,\gamma_1$.   Define the  realization set $\mcH_n=\mcH_{n,r}$ for $r \in \{2, \infty\}$  
\begin{align*}
    \mcH_n= \bigg \{ R_h(z,  \Lambda(z, p(z) )'\theta, p(z) ): \| p -p_0\|_r \leq p_{n}, \| \theta - \theta_0\|_1 \leq \theta_{n} \bigg \}.
\end{align*}

\item Suppose
 $$
R_I(z,\gamma_1)= \E   [\nabla_{\gamma_1} R(W,\gamma_1)|Z=z]
$$ is a  known function of $z,\gamma_1$.   Define the  realization set $I_n=I_{n,r}$ for $r \in \{2, \infty\}$ as 
\begin{align*}
    I_n=\bigg\{ R_I(w,p(z)): \| p -p_0\|_r \leq \pi_n \bigg \}
\end{align*}

\item  Suppose $m_{\text{pre}} (w,t,\gamma_1)$ is a known function of $w, t,\gamma_1$ and let $$q_0(z,t,\gamma_1) = \E[ m_{\text{pre}}(W,t,\gamma_1) | Z=z].$$ Define the realization set as
\begin{align*}
    \mcH_n(q_0)= \bigg \{ q(z, t, p(z)):  \sup_{t \in \R} \sup_{\gamma \in \Gamma} \| q(\cdot, t,\gamma)  - q_0(\cdot, t,\gamma) \|_r \leq q_{n}, \| \theta - \theta_0\|_1 \leq \theta_{n}, \| p - p_0 \|_r \leq \pi_n\bigg \}
\end{align*}

\end{enumerate}

  \begin{lemma}[Plausibility of Assumption \ref{ass:nuisance} for the adjusted CMR \eqref{eq:rho:phi}]
\label{lem:ver:nuisance}   Suppose Assumptions \ref{ass:convexitymoment}-\ref{ass:smooth} hold with $U \geq 1$. 
\begin{enumerate}
\item  The  plug-in estimator $\widehat{h}(z)$ 
\begin{align*}
\widehat{h}(z) := R_h(z, \Lambda(z, \widehat{p}(z))'\check{\theta}, \widehat{p}(z))
\end{align*}

of $h_0(z)$ belongs to $ \mcH_n$ w.p. $1-o(1)$.  The set $\mcH_n$ shrinks at rate $h_{n,r} = U^2 (\theta_n + \pi_{n,r})$ for $r= \infty$ and $h_{n,2} = \sqrt{2}U^2 (\theta_n + \pi_{n,2})$.
\item The plug-in estimator $\widehat{I}(z)$ 
\begin{align*}
 \widehat{I}(z)=R_I(z,\widehat{p}(z))
\end{align*}
of $I_0(z)$ belongs to  $I_n$ w.p. $1-o(1)$.  The set $I_n$ shrinks at rate $i_{n,r} = U  \pi_{n,r}$ for $r \in \{ \infty, 2 \}$.

\item  The estimator $\widehat{h}(z)$ of $h_0(z)$ 
\begin{align*}
\widehat{h}(z) := \widehat{q}(z, \Lambda(z,\widehat{p}(z))'\check{\theta},\widehat{p}(z))
\end{align*}
belongs to the realization set $ \mcH_n(q_0)=\mcH_{n,r}(q_0)$ for $r \in \{2, \infty\}$ w.p. $1-o(1)$.  The set $ \mcH_n(q_0)$ shrinks at rate $h_{n,r} = U^2 (\theta_n + \pi_{n,r})+ q_{n,r}$ for $r = \infty$ and $h_{n,2} =2 (U^2 (\theta_n + \pi_{n,2})+ q_{n,2})$. 

\end{enumerate}
\end{lemma}

\begin{proof}[Proof of Lemma \ref{lem:ver:nuisance}]

\textbf{ Step 1. } Proof of Lemma \ref{lem:ver:nuisance} (1).  By Assumption \ref{ass:smooth}, $m_{\text{pre}} (w,t, \gamma_1)$ and $\Lambda(z, \gamma_1)$ are smooth functions of $t$ and $\gamma_1$. Therefore, $R_h(z, t, \gamma_1)$ is a smooth function of $t$ and $\gamma_1$ whose partial derivatives are $U$-bounded uniformly over all arguments.  By intermediate value theorem, 
\begin{align*}
  R_h(z, t, \gamma) -  R_h(z, t_0, \gamma_0) &= \nabla_{t} R_h( z, \bar{t},\gamma) \Lambda(z, \gamma)'(\theta-\theta_0) + \nabla_{\gamma} R_h(z,\Lambda(z,\bar{\gamma})'\theta_0, \bar{\gamma}) \cdot (\gamma - \gamma_0) \\
  &= I_1(z) + I_2(z).
\end{align*}
By Cauchy Schwartz and Assumption \ref{ass:smooth}, for any $\theta: \| \theta - \theta_0\|_1 \leq \theta_{n}$,
\begin{align*}
\sup_{z \in \mcZ} | I_1(z) | &\leq U \sup_{z \in \mcZ} \sup_{\gamma \in \Gamma} | \Lambda(z, \gamma)'(\theta-\theta_0)  | \leq U^2 \| \theta - \theta_0 \|_1 \leq U^2 \theta_n .
\end{align*}
Plugging $\gamma = p(z)$ and $\gamma_0 = p_0(z)$ gives the bound in $\ell_{\infty}$ norm:
\begin{align*}
\sup_{z \in \mcZ} |I_1(z) + I_2(z) | \leq U^2 \theta_n + U \pi_n \leq U^2 (\theta_n + \pi_n).
\end{align*}
Therefore, $\mcH_n$ shrinks at rate $h_n = U^2 (\theta_n +  \pi_n)$ for $\pi_{n} = \pi_{n, \infty}$. \\
In $\ell_2$ norm, the bound is 
\begin{align*}
(\E (I_1(Z) + I_2(Z))^2)^{1/2} \leq \sqrt{2} ((\E I_1^2(Z))^{1/2}+ (\E I_2^2(Z))^{1/2}) \leq \sqrt{2} U^2 (\theta_n + \pi_n)
\end{align*}
for $\pi_{n} = \pi_{n, 2}$. 

\textbf{ Step 2. } Proof of Lemma \ref{lem:ver:nuisance} (2). By intermediate value theorem, $I_n$ in $\ell_r$ norm shrinks at rate $i_n = U  \pi_n$ if $\mcP_n$ shrinks at rate $\pi_{n,r}$, for $r \in \{2, \infty\}$. 

\textbf{ Step 3. } Proof of Lemma \ref{lem:ver:nuisance} (3). Observe that
\begin{align*}
 q(z, t, p(z)) -  q_0(z, t_0, p_0(z)) &=  (q(z, t, p(z)) -  q_0(z, t, p(z))) + (q_0(z, t, p(z))  -  q_0(z, t_0, p_0(z))) \\
 &= Q_1(z) + Q_2(z).
\end{align*}
Consider the case $r=\infty$. The bound on $Q_1(z)$ in $\ell_{\infty}$ follows from
 \begin{align*}
 \sup_{z \in \mcZ }  | Q_1(z) | \leq \sup_{t \in \R} \sup_{\gamma \in \Gamma} | q(z, t, \gamma) - q_0(z, t, \gamma) | \leq q_n.
 \end{align*}
Invoking Step 1 with $R_h(z, t, \gamma) = q_0(z,t,\gamma)$ gives
  \begin{align*}
 \sup_{z \in \mcZ }  | Q_2(z) | \leq U^2 (\theta_n + \pi_n).
 \end{align*}
 Consider the case $r=2$. The bound on $Q_1(z)$ in $\ell_2$ follows from
 \begin{align*}
(\E  Q_1^2(Z))^{1/2}  \leq \sup_{t \in \R} \sup_{\gamma \in \Gamma} (\E  (q(Z, t, \gamma) - q_0(Z, t, \gamma))^2)^{1/2} \leq q_n.
 \end{align*}
Invoking Step 1 with $R_h(z, t, \gamma) = q_0(z,t,\gamma)$ gives
  \begin{align*}
(\E  Q_2^2(Z))^{1/2}  \leq \sqrt{2} U^2 (\theta_n + \pi_n).
 \end{align*}

\end{proof}

\label{sec:ortho:proof}

\begin{proof}[Proof of Corollary \ref{cor:non-linear-te}]

\textbf{ Step 1. } Assumption \ref{ass:identification}  is a restatement of Assumption \ref{ass:reg:nonlinearte} \eqref{itm:nonlinear:ide}.
Assumption \ref{ass:identification2}  follows from Assumption \ref{ass:reg:nonlinearte} \eqref{itm:nonlinear:smooth}. Indeed, for any $w \in \mcW$ and $\nu \in {\cal C}(T;3)$ the following bound holds:
 \begin{align*}
 |(d - p_0(x)) \cdot (1,x)'(\theta_0 +  \nu)+q_0(x)| &\leq H_{\text{TE}}^2 \| \theta_0 +  \nu \|_1 + H_{\text{TE}} \leq 2 H_{\text{TE}}^3.
 \end{align*}
Therefore, 
\begin{align*}
  &\inf_{\theta \in \mathbb{B}}  \nabla_t \E[m(W,\Lambda(Z,g_0(Z))'\theta, g_0(Z))|Z=z] = \inf_{\theta \in \mathbb{B}} \frac{G' ((d - p_0(x)) \cdot (1,x)'\theta+q_0(x))}{V_0(z)} \geq  U_{\text{TE}} \text{B}_{\text{min}}.
\end{align*}
Assumption \ref{ass:convexitymoment} follows from Assumption  \ref{ass:reg:nonlinearte} \eqref{itm:nonlinear:mon}.

\textbf{ Step 2. } Verification of Assumption \ref{ass:smooth}. Let $x_0 = 1$.  For each $j \in \{0,1,2,\dots, p-1\}$  the function $\Lambda_j (z, \gamma_1) = (d - \gamma_1) \cdot x_j$ is bounded by $H_{\text{TE}} ^2$ for any $w \in \mcW$ and $\gamma_1 \in \Gamma_1$.  The first derivative $\nabla_{\gamma_1} \Lambda_j (z, \gamma_1) = x_j$ is bounded by $H_{\text{TE}}$ a.s. Finally, the second derivative $\nabla_{\gamma_1 \gamma_1} \Lambda_j (z, \gamma_1) = 0$ is zero. Next, the function $t  \rightarrow G(t) - y$ is $3$-times differentiable in $t$ with derivatives bounded by $U_{\text{TE}}$. The first (second) derivative of $\gamma_3^{-1}$ is bounded by   $U_{\text{TE}}$($2 U_{\text{TE}}$), respectively. Since the moment function $m(w, t, \gamma)$ is a product of $G(t + \gamma_2)-y$ and $\gamma_3^{-1}$, all the quantities in Assumption \ref{ass:smooth} are bounded by $U=C_{U_{\text{TE}}} U_{\text{TE}}^2$ for a sufficiently large absolute constant $C_{U_{\text{TE}}}$.

\textbf{ Step 3. } Verification of Assumption \ref{ass:nuisance}.   For $r \in \{2, \infty\}$, let  $\pi_n=\pi_{n,r}$ and $q_n=q_{n,r}$.  By Theorem \ref{cor:non-ortho-rate}, Step 1 of Definition \ref{def:twostep} converges in $\ell_1$-norm with $\theta_n = k C_{\text{TE}} \left(\sqrt{\frac{\log p}{n}} + \pi_n + q_n \right)$  for a sufficiently large $C_{\text{TE}}$.  By Lemma \ref{lem:ver:nuisance}(1), the estimate $\widehat{V}(d,x)$ in \eqref{eq:vee} converges at rate $ C_{U_{\text{TE}}} U_{\text{TE}}^2 ( \pi_n + q_n + \theta_n)$ in $\ell_r$ norm. Then, Assumption \ref{ass:nuisance} holds with $g_n = C_g k  \left( \pi_n + q_n+\sqrt{\frac{\log p}{n}} \right)$ for a sufficiently large absolute constant $C_g$.

\textbf{ Step 4. }  Verification of orthogonality condition \eqref{eq:orthogfinal}.   We show that the loss function gradient \eqref{eq:mom} is an orthogonal moment equation with respect to perturbations of $g_0(z)= \{ p_0(x), q_0(x), V_0(d,x)\}$. 
\begin{align*}
&\dfrac{\partial}{\partial r} \nabla_{\theta}  Q ( \theta_0, r (p - p_0) + p_0) = -\E  (D-p_0(X))   (p (X) -p_0(X)) \cdot (1,X)' \\
&- \dfrac{1}{V_0(D,X)} \E (Y - G ( (D-p_0(X)) \cdot ((1,X)'\theta_0) + q_0(X)) )  (p (X) -p_0(X)) \cdot (1,X)' =0 \\
&\dfrac{\partial}{\partial r} \nabla_{\theta}  Q ( \theta_0, r (q - q_0) + q_0) = -\E  (D-p_0(X))   (q (X) -q_0(X)) \cdot (1,X)' =0 \\
&\dfrac{\partial}{\partial r} \nabla_{\theta}  Q ( \theta_0, r (V - V_0) + V_0)  = -\dfrac{1}{V^2_0(D,X)} \E (Y - G ( (D-p_0(X)) \cdot ((1,X)'\theta_0) + q_0(X)) )\cdot \\
& (D-p_0(X)) (V(D,X) - V_0(D,X))\cdot (1,X)' =0.
\end{align*}

\end{proof}

\begin{proof}[Proof of Corollary \ref{cor:selection}]

The loss function \eqref{eqn:m-loss-md} in Example \ref{ex:MD}  is a special case of \eqref{eq:ortho}. By Theorem \ref{lem:general}, the population gradient based on \eqref{eqn:m-loss-md} obeys the orthogonality condition \eqref{eq:orthogfinal}.

\textbf{ Step 1.} Assumption \ref{ass:identification}  holds by Assumption \ref{ass:selection} \eqref{index-sel:identification}. Assumption \ref{ass:identification2} holds
by 
\begin{align*}
 \inf_{\theta \in \mathbb{B}}  \nabla_t \E[m(W,\Lambda(Z,g_0(Z))'\theta, g_0(Z))|Z=z] =\inf_{t \in \mathcal{T}} \nabla_t \E [u (Y, t)|X=x] \geq B_{\text{min}}
\end{align*}
for any $x \in \mathcal{X}$. Assumption \ref{ass:convexitymoment} directly follows from Assumption \ref{ass:selection} \eqref{index-sel:convexity}.

\textbf{ Step 2.} Verification of Assumption \ref{ass:smooth}.  Observe that $\Lambda(x, \gamma) = x$ and $\| x \|_{\infty} \leq H_{\text{MD}}$.  Next, the function $t  \rightarrow v \cdot u(y, t)$ is $3$-times differentiable in $t$ with derivatives bounded by $U_{\text{MD}}$. The second derivative of $\gamma_1^{-1}$ on $\Gamma_1$ is bounded by $6\bar{p}^{-3}$. Since the moment function $m_{\text{pre}}(w, t, \gamma_1)$ is a product of $v \cdot u(y,t)$ and $\gamma_1^{-1}$, all the quantities in Assumption \ref{ass:smooth} are bounded by $C_{\text{MD}} (U_{\text{MD}}^2 + 6\bar{p}^{-3})$ for a sufficiently large absolute constant $C_{\text{MD}}$.

\textbf{ Step 3.} Verification of Assumption \ref{ass:nuisance}.  For $r \in \{2, \infty\}$, let  $\pi_n=\pi_{n,r}$ and $q_n=q_{n,r}$. By Theorem \ref{cor:non-ortho-rate}, Step 1 of Definition \ref{def:twostep} converges in $\ell_1$-norm with $\theta_n = k C_{\text{MD}} \left(\sqrt{\frac{\log p}{n}} + \pi_n  \right)$ for a sufficiently large $C_{\text{MD}}$. By Lemma \ref{lem:ver:nuisance}(3), the estimate $\widehat{h}(x)$ in Remark \ref{rm:fsvee:selection}  converges at rate $C_{\text{MD}}  ( \pi_n  + \theta_n) + q_n$. Thus, Assumption \ref{ass:nuisance} holds with $$g_n = C_g \left(k  \left( \pi_n +\sqrt{\frac{\log p}{n}} \right)+q_n \right)$$ for a sufficiently large absolute constant $C_g$.

\end{proof}

\begin{proof}[Proof of Corollary \ref{cor:games}]

The loss function \eqref{eqn:m-loss-md} in Example \ref{ex:games}  is a special case of \eqref{eq:ortho}.  By Theorem \ref{lem:general}, the population gradient based on \eqref{eqn:m-loss-md} obeys the orthogonality condition \eqref{eq:orthogfinal}.

\textbf{ Step 1.} Assumption \ref{ass:identification}  holds by Assumption \ref{ass:games} \eqref{index-games:identification}. Assumption \ref{ass:identification2} holds
by 
\begin{align*}
 \inf_{\theta \in \mathcal{B}}  \nabla_t \E[m(W,\Lambda(Z,g_0(Z))'\theta, g_0(Z))|Z=z] = \inf_{t \in \mathcal{T}} \mcL (t)   \geq B_{\text{min}}
\end{align*}
for any $x \in \mathcal{X}$. Assumption \ref{ass:convexitymoment} follows from the monotonicity of logistic function.

\textbf{ Step 2.} Verification of Assumption \ref{ass:smooth}.  Observe that $\Lambda(x, \gamma) = (x;\gamma_1)$, $\| x \|_{\infty} \leq H_{\text{games}}$, and $\gamma_1 \in [0,1]$.  Next, the function $t  \rightarrow -(y - \mcL(t))$ is $3$-times differentiable in $t$ with derivatives bounded by some constant $U_{\text{games}}$, which can be taken $U_{\text{games}} \geq  H_{\text{games}}$. Therefore, all the quantities in Assumption \ref{ass:smooth} are bounded by $C_{U_{\text{games}}} U_{\text{games}}$ for a sufficiently large absolute constant $C_{U_{\text{games}}}$.

\textbf{ Step 3.} Verification of Assumption \ref{ass:nuisance}.  For  $r \in \{2, \infty\}$, let  $\pi_n=\pi_{n,r}$.  By Theorem \ref{cor:non-ortho-rate}, Step 1 of Definition \ref{def:twostep} converges in $\ell_1$-norm with   $\theta_n = k C_{\text{games}} \left(\sqrt{\frac{\log p}{n}} + \pi_n  \right)$ for a sufficiently large $C_{\text{games}}$.   By Lemma \ref{lem:ver:nuisance}(1), the estimate $\widehat{h}(x)$ in Remark \ref{rm:fsvee:games}  converges at rate $C_{\text{games}}  ( \pi_n  + \theta_n)$. Thus, Assumption \ref{ass:nuisance} holds with
$$g_n = C_g \left(k  \left( \pi_n +\sqrt{\frac{\log p}{n}} \right) \right)$$ for a sufficiently large absolute constant $C_g$.

\end{proof}

 \bibliography{my_new_bibtex}
 
\end{document}